\documentclass[12pt,reqno,final]{amsart}
\usepackage{amssymb,amsmath,amsthm,latexsym,graphicx}
\usepackage{latexsym,graphicx}
\usepackage[color]{showkeys}
\usepackage{appendix}
\definecolor{refkey}{gray}{.75}
\usepackage{fancyhdr}                %

\usepackage{hyperref}
\hypersetup{linktocpage}
\hypersetup{
    colorlinks,
    citecolor=black,
    filecolor=black,
    linkcolor=black,
    urlcolor=black
}

\newtheorem{theorem}{Theorem}[section]

\theoremstyle{remark}
\newtheorem{remark}{Remark}[section]

\numberwithin{equation}{section}
\numberwithin{figure}{section}

\def\clap#1{\hbox to 0pt{\hss#1\hss}}

\def\p{\partial}
\def\tilde{\widetilde}
\def\hat{\widehat}

\def\ud{\,\mathrm{d}}

\def\bk{\mathbf{k}}

\def\ub{\mathbf{u}}
\def\vb{\mathbf{v}}

\begin{document}
\title{Partial and spectral-viscosity models for geophysical flows}
\author{Qingshan Chen, Max Gunzburger and Xiaoming Wang}
\dedicatory{Dedicated to Roger Temam on the occasion of his seventieth birthday}
\thanks{This article appeared in {\em Chinese Annal of Mathematics
    Ser. B}, 31:579-606, 2010}
\thanks{DOI: 10.1007/s11401-010-0607-2}
\thanks{Corresponding author: Qingshan Chen ({\tt qchen3@fsu.edu})}

\date{\today}

\maketitle

\begin{abstract}

    Two models based on the hydrostatic primitive equations are 
    proposed. The first model is the primitive equations with 
    {\it partial} viscosity {\it only}, and is oriented towards
    large-scale wave structures in the ocean and atmosphere. The
    second model is the viscous primitive equations with spectral 
    eddy viscosity, and is oriented towards turbulent geophysical
    flows. For both models, the existence and uniqueness of 
    global strong solutions is established. For the second model,
    the convergence of the solutions to the solutions of the classical
    primitive equations as eddy viscosity parameters tend to zero 
    is also established.

\end{abstract}

\section{Introduction}\label{s1}
We study two models for geophysical flows based on the hydrostatic primitive equations; both are designed to faithfully simulate certain
phenomena in the geophysical flows but they are motivated
by different physical considerations.
A distinctive characteristic of the
flows under consideration is that
 the vertical scale ($\sim$ 10km)
is much smaller than the horizontal scale ($\sim$ 6000km).
Thanks to this
disparity, a hydrostatic
approximation is possible, and gives rise of the
primitive-equations.

On the mathematical side, the theory for the primitive
equations is fairly complete; see, e.g., the pioneering work in
\cite{LTW92a,LTW92b}, and the survey article
\cite{PTZ08}.
In particular, in contrast with
the Navier-Stokes equations \cite{ConFo88,Te01},
the primitive equations have been
shown to have unique global strong solutions \cite{CT07, Ko07,
KuZi07, Ju07}.

The first model we study aims to faithfully
simulate large-scale coherent structures including wave phenomena
in the ocean and atmosphere. For relevant discussions of this topic,
see, e.g., \cite{gill,MaWa06,P87}.
 For the phenomena that we are interested in, both the ocean and
atmosphere are close to being inviscid. Therefore, the inviscid
primitive equation is the preferred model. However, it is very
costly to simulate the inviscid primitive equations directly due to
the small scales embodied within the model. What we propose here is
a new model with eddy viscosity added to the small scale (high
frequency) part only while keeping the large scale (low frequency)
features intact (at least directly). Intuitively, this type of model
would reduce the complexity due to the damping on small
scales whereas keeping the desirable large scale structures. Such a naive
approach may not work all the time due to the cascade of energy
induced by the nonlinear advection term. However, we can easily find
situations where such cascade is small or negligible. Indeed, it is
easy to find exact large-scale solutions to this primitive equation
with partial viscosity; see below. Of course, the existence of such
examples do not fully justify the model, and extensive numerical
experiments are called for which is our future plan. The idea of
partial damping on the high frequency components of the system is
not alien to the geophysical community (see, e.g., \cite{FDK96,
Mcw84}) or the mathematical community (see e.g. ~\cite{NeTa92,
CLG05, CLG07, DGK08}), although the application to the primitive
equations is new here. Our goal in this paper is to demonstrate the
global well posedness of this model with partial viscosity.

The second model
is oriented at turbulence modeling for geophysical flows.
The simulation of three-dimensional turbulent flows is a formidable task due to
the need to resolve the small scale fluctuations or {\em eddies}
that have subtle effects on the large-scale dynamics of the flow.
To make this problem computationally tractable, these effects
must be modeled. In one approach, the velocity field is averaged
over a small radius to derive equations in terms of the averaged
velocity. For nonlinear equations, there arises the
problem of {\em closure} because the product operation is not
closed under the averaging process. To obtain a closed system,
the average of the nonlinear terms in the equations must
be approximated and expressed solely in terms of averaged
quantities. The way in which this is done gives rise to a variety
of models. The approach we consider, called the {\it eddy-viscosity
method}, treats the Reynolds stress as a viscous effect caused
by the transport and dissipation of energy due to the small-scale
eddies. For this reason, this additional viscosity is called the
{\it eddy viscosity} or {\it turbulent viscosity}. The turbulence
model of Smagorinsky \cite{Sm63} belongs to this type. For an overall
survey on issues related to these models, see \cite{BeIlLa06}.

Unfortunately, a straightforward application of the approach
described above leads to the over smearing of the large-scale
structures in the fluid. To remedy this unwanted effect, it has
been proposed that the eddy viscosity be added only to the
{\it subgrid scales}. In this way, one hopes to prevent the
large-scale structure from being smeared away. Here, we examine
a particular class of models of this type called {\it
spectral-viscosity} or {\it spectral-vanishing-viscosity} models,
in which the scales are defined in terms
of Fourier modes. The subgrid viscosity is simply realized as an
addition of the artificial viscosity only to the high-frequency
modes. The most intuitive way of doing this to insert a high-pass
filter to the standard eddy viscosity.
This approach was considered in \cite{GuPr03} for hyperviscosity
on the Navier-Stokes equations and in \cite{GLSTW} for nonlinear
as well as hyper-viscosity on the Navier
Stokes Equations. In both works, the well-posedness of the resulting
spectral viscosity is proven, and the consistency of these
model with the original Navier-Stokes equations is discussed.

We employ the idea of {\it spectral viscosity}
to build and analyze a turbulence model for the geophysical flows
in the ocean and atmosphere.
As mentioned above, the primitive equations, even without
any eddy viscosity, have been shown to have unique global
solutions, provided the initial and boundary data are sufficiently
smooth. For our model, we prove its global well posedness,
which should not come as a surprise. In addition, we will show
the convergence of the solutions of the model
to the solutions of the primitive equations
without any eddy viscosity, as the eddy viscosity parameters
tend to zero.
This is not
possible yet for the Navier-Stokes equations because there
convergence is only shown to be in a weak sense.

We should point out that this technique, usually under the name of
{\it spectral viscosity} or {\it spectral-vanishing viscosity}, is
known in terms of turbulence modelling (see \cite{KaraKarn00,
GuPr03, AvX09, GLSTW,SSK05}, and below), and to applications in
geophysical fluid dynamics \cite{FDK96}. However, to the best of our
knowledge, the well-posedness result for the three-dimensional nonlinear primitive
equations with partial viscosity is new.

The paper is organized as follows.
In Section \ref{s3}, we introduce and prove
the well posedness of a model with only partial viscosity.
In Section \ref{s2}, we
introduce the linear spectral eddy-viscosity model. In Section
\ref{s2.2}, we prove the existence and uniqueness of strong solutions
to the model. In Section \ref{s2.3} we study the convergence of
the solutions of the model to the solutions of the original
primitive equations.

\section{A model with partial viscosity }\label{s3}
In this section we study a model with partial high--frequency
viscosity only. The lower modes are not damped directly. This
feature renders the model suitable for large scale coherent
structures in the ocean and atmosphere, because non--physical large
scale damping could change the large scale coherent structures over
time.

The model reads
\begin{align}
   &\dfrac{\p\ub}{\p t} + (\ub\cdot\nabla)\ub + w\dfrac{\p\ub}{\p z}
   +f\bk\times\ub + \dfrac{1}{\rho_0}\nabla p - { } \nonumber \\
   &\phantom{\dfrac{\p\ub}{\p t} + { }} \mu\Delta (I-P_{M,N})\ub
   - \nu\dfrac{\p^2}{\p z^2}(I-P_{M,N})\ub  = \mathbf{F},\label{e4.1}\\
  & \dfrac{\p p}{\p z} = -\rho_0g,\label{e4.2}\\
  & \nabla\cdot\ub + \dfrac{\p w}{\p z} = 0,\label{e4.3}
\end{align}
In the above, $\ub=(u,v)$ is the horizontal velocity, $w$ is the
vertical velocity, $\mu$ and $\nu$ are the horizontal and vertical
kinetic viscosities, respectively.
We
use $\nabla$ and $\Delta$ to denote the 2D horizontal gradient
and Laplacian operators, respectively. The operator $(I-P_{M,N})$
represents the high--pass filter and will be defined later on.

We consider a rectangular domain $\Omega = M\times(-H,0)$, with
$M=(0,L_x)\times(0,L_y)$. We consider periodic boundary conditions
in both the $x$ and $y$ directions, and free--slip, non--penetration
boundary conditions in the vertical directions. More precisely,
\begin{align}
   \ub(x,y,z,t) = \ub(x+L_x,y,z,t),\label{e4.3a}\\
   \ub(x,y,z,t) = \ub(x,y+L_y,z,t),\label{e4.3b}\\
   \dfrac{\p\ub}{\p n}|_{z=0} =
   \dfrac{\p\ub}{\p n}|_{z=-H} = 0,\label{e4.3c}\\
   w|_{z=0} = w|_{z=-H} = 0,\label{e4.3d}\\
   p|_{z=0} = p_0.\label{e4.3e}
\end{align}

Under the settings just described, we can define the high--pass
filter in terms of Fourier frequencies. Specifically,
for each function $\ub\in(L^2(\Omega))^2$ we let
\begin{equation}
   P_{M,N}\ub = \sum_{|\mathbf{m}|_{sup}\leq M\,\&\, n\leq N} \hat\ub_{\mathbf{m},n}
   e^{i\mathbf{m}\cdot\mathbf{x}'}\cos nz'
   \label{e4.3f}
\end{equation}
where
\begin{align*}
   &\mathbf{m} = (m_1,m_2)\in \mathbb{Z}^2,\\
   &|\mathbf{m}|_{\sup} = \max(|m_1|, |m_2|),\\
   &\mathbf{x}' = 2\pi(\dfrac{x}{L_x},\dfrac{y}{L_y}),\\
   &z' = \dfrac{\pi z}{H},\\
   &\hat{\ub}_{\mathbf{m},n} = \int_{\Omega} \ub(x,y,z)
   e^{-i\mathbf{m}\cdot\mathbf{x}'}\cos nz'\ud x\ud y\ud z.
\end{align*}

\begin{remark}
   The issue of suitable physical boundary conditions for the
   inviscid primitive equations is an unresolved one. Partial results concerning
   the linearized primitive equations
   are available in a series of papers \cite{RTT05a}, \cite{CLRTT07}, and
   \cite{RTT08}. There, an infinite set of nonlocal boundary conditions
   were proposed, which guaranteed the well--posedness of the linearized
   system. Here we avoid this issue and use the periodic boundary conditions
   on the lateral boundaries.
\end{remark}

\begin{remark}
    As we have touched upon in Introduction, we can easily identify
    some large--scale motions which are exact solutions of the inviscid
    primitive equations. For example, let
    \begin{displaymath}
    \psi = \sin\left(\dfrac{4\pi x}{L_1}\right)
    \cos\left(\dfrac{2\pi y}{L_2}\right),
    \end{displaymath}
    and
    \begin{displaymath}
    \ub = \nabla^{\perp}\psi.
    \end{displaymath}
    With the surface pressure $p_0$ given by
    \begin{displaymath}
    p_0 = f\psi + \dfrac{2\pi^2}{L_2^2}\cos\left(\dfrac{4\pi x}{L_1}\right)^2 -
    \dfrac{8\pi^2}{L_1^2}\cos\left(\dfrac{2\pi y}{L_2}\right)^2,
    \end{displaymath}
    $(\ub, 0)$ is a set of exact solutions of the inviscid primitive
    equations, that is, the system \eqref{e4.1}--\eqref{e4.3e} without
    the viscosities, or with partial viscosity with $M>4$.

\end{remark}

\subsection{The barotropic and baroclinic modes}\label{s3.1}

As usual, $p$ and $w$ can be expressed in terms of $\ub$. Specifically,
integrating \eqref{e4.2} from $z$ to 0, and using the boundary condition
\eqref{e4.3e} we obtain
\begin{equation}
   p(x,y,z,t) = p_0(x,y,t) - \rho_0gz,\label{e4.3g}
\end{equation}
Integrating \eqref{e4.3} from $z$ to 0, and using the boundary condition
\eqref{e4.3d}, we obtain
\begin{equation}
   w(x,y,z,t) = \int_z^0\nabla\cdot\ub\ud\xi = \nabla\cdot\int_z^0\ub\ud\xi.
   \label{e4.3h}
\end{equation}
Setting $z=-H$ in \eqref{e4.3h}, and using \eqref{e4.3d} again, we find
\begin{equation}
   \nabla\cdot\int_{-H}^0\ub\ud z = 0.
   \label{e4.3l}
\end{equation}
We substitute \eqref{e4.3g} and \eqref{e4.3h} into \eqref{e4.1}, and
obtain a single closed equation for $\ub$:
\begin{multline}
   \dfrac{\p\ub}{\p t} + (\ub\cdot\nabla)\ub +
   (\nabla\cdot\int_z^0\ub\ud\xi)\dfrac{\p\ub}{\p z}
   +f\bk\times\ub + \dfrac{1}{\rho_0}\nabla p_0 - \\
   \mu\Delta (I-P_{M,N})\ub
   - \nu\dfrac{\p^2}{\p z^2}(I-P_{M,N})\ub = \mathbf{F}.
   \label{e4.4}
\end{multline}
It turns out essential to rewrite this equation
and work with the following form:
\begin{multline}
   \dfrac{\p\ub}{\p t} + (\ub\cdot\nabla)\ub +
   \left(\nabla\cdot\int_z^0\ub\ud\xi\right)\dfrac{\p\ub}{\p z}
   +f\bk\times\ub + \dfrac{1}{\rho_0}\nabla p_0 - \\
   \mu\Delta \ub
   - \nu\dfrac{\p^2}{\p z^2}\ub = \mathbf{F} -
   \mu\Delta P_{M,N} \ub -
   \nu\dfrac{\p^2}{\p z^2} P_{M,N} \ub.
   \label{e4.4a0}
\end{multline}
We note here that $ P_{M,N} \ub$ contains
only finite number modes of $\ub$.
The prognostic variable $\ub$ satisfies the equation \eqref{e4.4a0},
the constraint \eqref{e4.3l}, and
the boundary conditions \eqref{e4.3a}--\eqref{e4.3c}. To complete the
system, we also require $\ub$ to satisfy the following initial
condition:
\begin{equation}
   \ub(x,y,z,0) = \ub_0(x,y,z).
   \label{e4.4a1}
\end{equation}
The diagnostic variables $w$ and $p$ are given by \eqref{e4.3h} and
\eqref{e4.3g} respectively.

We now specify the barotropic and baroclinic modes of equation
\eqref{e4.4a0}.
We let
\begin{equation}
   \overline{\ub}(x,y,t) = \dfrac{1}{H}\int_{-H}^0\ub(x,y,z,t)\ud z,
   \label{e4.4a2}
\end{equation}
which denotes the barotropic mode of the primitive variables. We
also denote by
\begin{equation}
   \ub'(x,y,z,t) = \ub(x,y,z,t)-\overline{\ub}(x,y,t)
   \label{e4.4a_2}
\end{equation}
the baroclinic mode of the primitive variables. It is easy to see that
\begin{equation}
   \overline{\ub'} = 0.
   \label{e4.4a3}
\end{equation}
From \eqref{e4.3a} and \eqref{e4.3b} we derive the boundary conditions for
$\overline{\ub}$:
\begin{align}
   \overline{\ub}(x,y,t) = \overline{\ub}(x+L_x,y,t),\label{e4.4a4}\\
   \overline{\ub}(x,y,t) = \overline{\ub}(x,y+L_y,t).\label{e4.4a5}
\end{align}
By \eqref{e4.3l} we see that $\overline{\ub}$ satisfies the following
constraint:
\begin{equation}
   \nabla\cdot\overline{\ub} = 0.
   \label{e4.4a6}
\end{equation}
Then $\ub'$ also satisfy the periodic boundary conditions on the
lateral boundary,
\begin{align}
   \ub'(x,y,z,t) = \ub'(x+L_x,y,z,t),\label{e4.4a7}\\
   \ub'(x,y,z,t) = \ub'(x,y+L_y,z,t).\label{e4.4a8}
\end{align}
It inherits the boundary conditions for $\ub$ on the top
and bottom,
\begin{equation}
   \left.\dfrac{\p\ub'}{\p z}\right|_{z=0} = \left.\dfrac{\p \ub'}{\p z}\right|_{z=-H} = 0.
   \label{e4.4a9}
\end{equation}
We now derive the equations that $\overline{\ub}$ and $\ub'$ satisfy by first
taking average of the equation \eqref{e4.4a0}:
\begin{multline}
   \dfrac{\p\overline{\ub}}{\p t} + \overline{(\ub\cdot\nabla)\ub} +
   \overline{w\dfrac{\p\ub}{\p z}}
   +f\bk\times\overline{\ub} + \dfrac{1}{\rho_0}\nabla p_0 -
   \mu\Delta\overline{\ub}
   - \nu\overline{\dfrac{\p^2{\ub}}{\p z^2}}  = \\
   \overline{\mathbf{F}}
   -\mu\overline{\Delta P_{M,N}\ub} -
   \nu\overline{\dfrac{\p^2}{\p z^2}P_{M,N}\ub}.
   \label{e4.4a10}
\end{multline}
We notice that
\begin{align*}
 &\overline{(\dfrac{\p^2{\ub}}{\p z^2})} =
 \dfrac{1}{H}\int_{-H}^0 \dfrac{\p^2{\ub}}{\p z^2} \ud z =
 \dfrac{1}{H}\left.\dfrac{\p\ub}{\p z}\right|^0_{-H} = 0,
 \quad(\textrm{ by \eqref{e4.3c}}),\\
 &\overline{\dfrac{\p^2}{\p z^2}P_{M,N}\ub}
   = \dfrac{1}{H}\int_{-H}^0
   \dfrac{\p^2}{\p z^2}P_{M,N}\ub\ud z
   =0,\\
   & \overline{\Delta P_{M,N}\ub} = \Delta P_{M,N}\overline{\ub}.
\end{align*}
By using \eqref{e4.4a3}, we find that
\begin{equation*}
   \overline{(\ub\cdot\nabla)\ub} =
   (\overline{\ub}\cdot\nabla)\overline{\ub} +
   \overline{(\ub'\cdot\nabla)\ub'}.
\end{equation*}
Using \eqref{e4.3l}, \eqref{e4.4a3} and \eqref{e4.4a6}, we find
\begin{equation*}
   \overline{w\dfrac{\p\ub}{\p z}} =
   \overline{(\nabla\cdot\ub')\cdot\ub'}.
\end{equation*}
Hence the equation for $\overline{\ub}$ can be written as
\begin{multline}
   \dfrac{\p\overline{\ub}}{\p t} +
   (\overline{\ub}\cdot\nabla)\overline{\ub} +
   [\overline{(\ub'\cdot\nabla)\ub'} +
   \overline{(\nabla\cdot\ub')\ub'} ] \\
   +f\bk\times\overline{\ub} +
   \dfrac{1}{\rho_0}\nabla p_0 -
   \mu\Delta\overline\ub =
   \overline{\mathbf{F}} - \mu\Delta P_{M,N} \overline{\ub}.\label{e4.4a}
\end{multline}
The barotropic variable $\overline{\ub}$ satisfies the following
conditions:
\begin{gather}
   \overline{\ub}(x,y,t) = \overline{\ub}(x+L_x,y,t),\label{e4.5}\\
   \overline{\ub}(x,y,t) = \overline{\ub}(x,y+L_y,t).\label{e4.6}\\
   \nabla\cdot\overline{\ub} = 0. \label{e4.7}
\end{gather}
Subtracting \eqref{e4.4a} from \eqref{e4.4a0} we obtain
the equation for the baroclinic mode:
\begin{multline}
   \dfrac{\p\ub'}{\p t} +
   (\ub'\cdot\nabla)\ub' +
   \left(\nabla\cdot\int_z^0\ub'\ud\xi\right)\dfrac{\p\ub'}{\p z}
   + f\bk\times\ub' - \mu\Delta\ub' - \nu\dfrac{\p^2\ub'}{\p z^2}\\
   { }+\left[(\ub'\cdot\nabla)\overline{\ub} + (\overline{\ub}\cdot\nabla)\ub'
   - \left(\overline{(\ub'\cdot\nabla)\ub'}+
   \overline{(\nabla\cdot\ub')\ub'}\right)\right]
   = \mathbf{F}'- \\
   \mu\Delta P_{M,N} \ub' -
   \nu\dfrac{\p^2}{\p z^2} P_{M,N} \ub'.
   \label{e4.8}
\end{multline}
In addition $\ub'$ satisfy
\begin{align}
   &\ub'(x,y,t) = \ub'(x+L_x,y,t),\label{e4.9}\\
   &\ub'(x,y,t) = \ub'(x,y+L_y,t),\label{e4.10}\\
   &\left.\dfrac{\p\ub'}{\p z}\right|_{z=0} = \left.\dfrac{\p \ub'}{\p z}\right|_{z=-H} = 0,
   \label{e4.11}\\
   &\int_{-H}^0 \ub'\ud z = 0.\label{e4.12}
\end{align}

\subsection{Global well--posedness of the model with partial viscosity}
In order to have uniqueness for the solution of \eqref{e4.1}--\eqref{e4.3e}
we work with functions that have zero average over $\Omega$. Indeed
it can be checked that if the initial data and the forcing have
zero average over $\Omega$, then the solutions have zero average
over $\Omega$ at any time.

We use the convention $\dot{L}^2(\Omega)$, $\dot{H}^1(\Omega)$,
etc.~to denote function spaces that have zero average over $\Omega$.
We let
\begin{align*}
   & H = (\dot{L}^2(\Omega))^2,\\
   &\mathcal{V} = \big\{\ub\in\dot{\mathcal{C}}^\infty(R^3)^2 |\,
   \ub \textrm{ periodic in $x$ with period $L_x$, periodic in $y$ with} \\
   &\phantom{\mathcal{V} = }\textrm{ period $L_y$, periodic and even in $z$ ith period $2H$}
   \big\},\\
   & V = \overline{\mathcal{V}}^{H^1} (\textrm{ closure of } \mathcal{V}
     \textrm{ in }(H^1(\Omega))^2 ).
\end{align*}
The inner product and norm of $H$ will be denoted as
$(\cdot,\,\cdot)$ and $|\cdot|$, respectively. The space
$V$ inherits the inner product and norm of $H^1$, which will be
denoted as $( (\cdot,\,\cdot) )$ and $||\cdot||$, respectively.

Since the functions in $V$ has zero spatial averages, we have
the following Poincar\'e inequality for functions in $V$:
\begin{equation}
   |\ub|^2\leq C(|\nabla\ub|^2 + |\dfrac{\p\ub}{\p z}|^2).
   \label{e34}
\end{equation}
Therefore $(|\nabla\ub|^2 + |\dfrac{\p\ub}{\p z}|^2)^{\frac{1}{2}}$
is equivalent to the usual $H^1$ norm, and can be taken as the
norm for $V$.

In what follows we abuse the notation by
denoting every generic constant by $C$. Such constants may depend on
the domain $\Omega$ and the function spaces in the context, but
we omit such dependence in the notation.
But if the constant depends on any other parameters, such as
$M$, $N$ etc.,
we shall use a specific symbol and specify such dependence
in the notation.

We shall  prove the global
existence and uniqueness of strong solutions to the
system \eqref{e4.1}--\eqref{e4.3e}.
\begin{theorem}\label{t0}
   For a given $T>0$, let
   $\mathbf{F}\in L^2(0,T;H)$,
   $\ub_0\in V$. Then there exists a unique strong solution
   $\ub\in \mathcal{C}([0,T];V)\cap L^2(0,T; H^2(\Omega))$ of the
   system \eqref{e4.1}--\eqref{e4.3e} which depends continuously
   on the initial data.
\end{theorem}

We shall first obtain some key estimates that will be needed
for the proof of Theorem \ref{t0}. The proof of the theorem will
be furnished at the end of this subsection.

\vspace{5mm}
\noindent{\it $L^2$ estimates}\\
We multiply \eqref{e4.4a0} by $\ub$, integrate by parts
over $\Omega$, using the boundary conditions
\eqref{e4.3a}--\eqref{e4.3e}, we obtain
\begin{multline}
   (\dfrac{\p\ub}{\p t},\ub) + \left( (\ub\cdot\nabla)\ub,\ub \right) +
   \left( (\nabla\cdot\int_z^0\ub\ud\xi)\dfrac{\p\ub}{\p z},\,\ub \right) + { }\\
   (f\bk\times\ub,\ub) + \dfrac{1}{\rho_0}(\nabla p_0,\ub) +
   \mu|\nabla\ub|^2 + \nu\left|\dfrac{\p \ub}{\p z}\right|^2  \\
   =
   (\mathbf{F},\ub) +
   \mu_\delta|\nabla P_{M,N}\ub|^2 +
   \nu_\delta\left|\dfrac{\p}{\p z}P_{M,N}\ub\right|^2.
   \label{e4.12a}
\end{multline}
We notice that
\begin{align}
   &\left(\dfrac{\p\ub}{\p t},\,\ub\right) =
   \dfrac{1}{2}\dfrac{\ud}{\ud t}|\ub|^2,\label{e4.12b}\\
   &(\nabla p_0,\,\ub)_\Omega = H(\nabla p_0,\,\overline{\ub})
   = 0,\label{e4.12c}\\
   &(f\bk\times\ub,\ub) = 0.\label{e4.12d}
\end{align}
Let
\begin{equation}
   b(\ub,\tilde\ub,\ub^\#) =
   \left( (\ub\cdot\nabla\tilde\ub) +
   w(\ub)\dfrac{\p\tilde\ub}{\p z},\,\ub^\# \right),
   \label{e4.12e}
\end{equation}
where $w(\ub)$ is defined as in \eqref{e4.3h}. We can verify that
the trilinear operator $b(\ub,\tilde\ub,\ub^\#)$ is skew symmetric
with respect to the last two arguments,
that is,
\begin{equation}
   b(\ub,\tilde\ub,\ub^\#) = -b(\ub,\ub^\#,\tilde\ub).
   \label{e4.12f}
\end{equation}
Then it is inferred from \eqref{e4.12f} that
\begin{equation}
   b(\ub,\ub,\ub) = 0.
   \label{e4.12g}
\end{equation}
By \eqref{e4.12a}, \eqref{e4.12b}, \eqref{e4.12c},
\eqref{e4.12d} and \eqref{e4.12g} we have
\begin{multline}
   \dfrac{1}{2}\dfrac{\ud}{\ud t}|\ub|^2
   + \mu |\nabla \ub|^2 +
   \nu |\dfrac{\p}{\p z}\ub|^2 =
   (\mathbf{F},\ub)+\\
   \mu |\nabla P_{M,N} \ub|^2 +
   \nu |\dfrac{\p}{\p z} P_{M,N} \ub|^2.
   \label{e4.13}
\end{multline}
Since $ P_{M,N} \ub$ contains only finite number of
Fourier modes,$|\nabla P_{M,N} \ub|^2$, as well as
$|\dfrac{\p}{\p z} P_{M,N} \ub|^2$ can be bounded
by $|\ub|^2$. More generally, for functions in $H^k(\Omega)$,
there exists a constant $C_k(M,N,\Omega)$, which is independent
of the function, such that
\begin{equation}
   | P_{M,N} \ub|_{H^k} \leq C_k |\ub|_{L^2}.
   \label{e4.14}
\end{equation}
Using \eqref{e4.14} and the Cauchy--Schwarz inequality,
we derive from \eqref{e4.13} that
\begin{equation}
   \dfrac{\ud}{\ud t}|\ub|^2
   + \mu |\nabla \ub|^2 +
   \nu |\dfrac{\p}{\p z}\ub|^2\leq
   C|\mathbf{F}|^2 + C_1|\ub|^2.
   \label{e4.15}
\end{equation}
Applying the Gronwall inequality to \eqref{e4.15}
yields
\begin{equation}
   |\ub(\cdot,t)|^2 + \mu\int_0^t|\nabla\ub|^2\ud s +
   \nu\int_0^t\left|\dfrac{\p}{\p z}\ub\right|^2\ud s
    \leq J_1(t),
   \label{e4.16}
\end{equation}
where
\begin{equation}
   J_1(t) \equiv e^{C_1 t}\left(
   |u_0|^2 + C\int_0^t|\mathbf{F}|^2\ud s\right).
   \label{e4.17}
\end{equation}

\vspace{5mm}
\noindent{\it $L^6$ estimate on $\ub'$}\\
We take inner product of \eqref{e4.8} with $|\ub'|^4\ub'$,
and integrate by parts over $\Omega$ to obtain
\begin{multline}
   \dfrac{1}{6}\dfrac{\ud}{\ud t}|\ub'|_{L^6}^6 +
   b(\ub',\ub',|\ub'|^4\ub') +
   \int_\Omega f\bk\times\ub'\cdot|\ub'|^4\ub'\ud\Omega + { }\\
   \int_\Omega\left( (\ub'\cdot\nabla)\overline\ub +
   (\overline\ub\cdot\nabla)\ub' -
   \overline{(\ub'\cdot\nabla)\ub' +
   (\nabla\cdot\ub')\ub'} \right)
   |\ub'|^4\ub'\ud\Omega - { }\\
   \int_\Omega\left(\mu\Delta\ub' +
   \nu\dfrac{\p^2\ub'}{\p z^2}\right)|\ub'|^4\ub'\ud\Omega 
   = \int_\Omega\mathbf{F}'\cdot|\ub'|^4\ub'\ud\Omega - { }\\
      \int_\Omega\left(\mu\Delta P_{M,N} \ub' +
      \nu\dfrac{\p^2}{\p z^2} P_{M,N} \ub'\right)|\ub'|^4\ub'\ud\Omega
   \label{e4.18}
\end{multline}
In the above, the trilinear operator $b(\cdot,\,\cdot,\,\cdot)$ is
defined as in \eqref{e4.12e}.
We can verify by calculations
that
\begin{equation}
   b(\ub',\ub',|\ub'|^4\ub') = -2 b(\ub',\ub',|\ub'|^4\ub'),
   \label{e4.18a}
\end{equation}
which implies that
\begin{equation}
   b(\ub',\ub',|\ub'|^4\ub') = 0.
   \label{e4.18b}
\end{equation}
Since $k\times\ub'$ is orthogonal to $\ub'$, we have
\begin{equation}
 \int_\Omega f\bk\times\ub'\cdot|\ub'|^4\ub'\ud\Omega = 0.
   \label{e4.18c}
\end{equation}
Noticing the divergence free condition \eqref{e4.7} for
$\overline{\ub}$ and the horizontal periodic boundary conditions
\eqref{e4.5}--\eqref{e4.6} and \eqref{e4.9}--\eqref{e4.10} for
$\overline{\ub}$ and $\ub'$ respectively, we find
\begin{equation}
   \int_\Omega(\overline\ub\cdot\nabla)\ub'\cdot
   |\ub'|^4\ub'\ud\Omega = 0.
   \label{e4.18d}
\end{equation}
For the inner products involving the diffusion terms,
we find
\begin{equation}
   -\int_\Omega\mu\Delta\ub'\cdot|\ub'|^4\ub'\ud\Omega =
   \mu\int_\Omega |\nabla\ub'|^2|\ub'|^4\ud\Omega +
   \mu\int_\Omega\left|\nabla|\ub'|^2\right|^2|\ub'|^2\ud\Omega,
   \label{e4.18e}
\end{equation}
\begin{equation}
   -\int_\Omega\nu\dfrac{\p^2\ub'}{\p z^2}\cdot|\ub'|^4\ub'\ud\Omega =
   \nu\int_\Omega |\dfrac{\p \ub'}{\p z}|^2|\ub'|^4\ud\Omega +
   \nu\int_\Omega\left|\dfrac{\p}{\p z}|\ub'|^2\right|^2|\ub'|^2\ud\Omega.
   \label{e4.18f}
\end{equation}
For integrals on the right--hand
side of \eqref{e4.18}, we use \eqref{e4.14} (with $k=2$)
to find that
\begin{align*}
   &\phantom{\leq}\int_\Omega -\mu\Delta P_{M,N} \ub'\cdot
   |\ub'|^4\ub'\ud\Omega\\
   &\leq \mu|\Delta P_{M,N} \ub'|_{L^2}
   |\ub'|_{L^{10}}^5\\
   &\leq C_2|\ub'|_{L^2}\left| |\ub'|^3\right|_{L^{\frac{10}{3}}}^
   {\frac{5}{3}}\\
   &\leq C|\ub'|_{L^2}^2\left| \ub'\right|_{L^6}^{4} +
    \dfrac{\mu}{4}\int_\Omega|\ub'|^4|\nabla\ub'|^2\ud\Omega
    + \dfrac{\nu}{4}\int_\Omega|\ub'|^4\left|\dfrac{\p\ub'}{\p z}\right|^2\ud\Omega.
\end{align*}
In the above we have used
 the
interpolation inequality
\begin{equation}
   |\phi|_{L^\frac{10}{3}} \leq C|\phi|_{L^2}^{\frac{2}{5}}
   |\phi|_{H^1}^{\frac{3}{5}},
   \label{e4.18g}
\end{equation}
which can be obtained by setting $p = 10/3,\,p_1 = 2,\,p_2 = 6$
in \eqref{a1}, and then using \eqref{a5}.
Similarly,
\begin{multline*}
   \int_\Omega -\nu\dfrac{\p^2}{\p z^2} P_{M,N} \ub'\cdot
   |\ub'|^4\ub'\ud\Omega\\
   \leq C|\ub'|_{L^2}^2\left| \ub'\right|_{L^6}^{4} +
    \dfrac{\mu}{4}\int_\Omega|\ub'|^4|\nabla\ub'|^2\ud\Omega
    + \dfrac{\nu}{4}\int_\Omega|\ub'|^4\left|\dfrac{\p\ub'}{\p z}\right|^2\ud\Omega.
\end{multline*}
Hence we derive from \eqref{e4.18} that
\begin{multline}
    \dfrac{1}{6}\dfrac{\ud}{\ud t}|\ub'|_{L^6}^6 +{ }\\
   \int_\Omega\left( (\ub'\cdot\nabla)\overline\ub -
   \overline{(\ub'\cdot\nabla)\ub' +
   (\nabla\cdot\ub')\ub'} \right)
   \cdot|\ub'|^4\ub'\ud\Omega + \\
   \dfrac{1}{2}\mu\int_\Omega\left(|\nabla\ub'|^2|\ub'|^4 +
   \left|\nabla|\ub'|^2\right|^2|\ub'|^2\right)\ud\Omega + { }\\
   \dfrac{1}{2}\nu\int_\Omega\left(\left|\dfrac{\p \ub'}{\p z}\right|^2|\ub'|^4 +
   \left|\dfrac{\p}{\p z}|\ub'|^2\right|^2|\ub'|^2 \right)\ud\Omega +\\
   \leq \int_\Omega\mathbf{F}'\cdot|\ub'|^4\ub'\ud\Omega +
   C |\ub'|_{L^2}^2|\ub'|^4_{L^6}.
   \label{e4.19}
\end{multline}
For the integrals on the right--hand side of
\eqref{e4.19} that involve nonlinear terms, we proceed
by integration by parts, using the periodic boundary conditions on $\ub'$ and
$\overline{\ub}$ when appropriate, and we find that
\begin{align*}
   &\phantom{=}\left|\int_{\Omega}(\ub'\cdot\nabla)\overline{\ub}\cdot|\ub'|^4\ub'\ud\Omega\right|\\
   &= \left| \int_{\Omega}\left( (\nabla\cdot\ub')\overline{\ub}\cdot|\ub'|^4\ub' +
     (\ub'\cdot\nabla)(|\ub'|^4\ub')\cdot\overline{\ub}\right)\ud\Omega\right|\\
    &\leq \left|\int_{M}\left(\overline{\ub}\cdot
       \int_{-H}^0 (\nabla\cdot\ub')|\ub'|^4\ub'\ud z +
       \overline{\ub}\cdot
       \int_{-H}^0(\ub'\cdot\nabla)(|\ub'|^4\ub')\ud z\right)
       \ud x\ud y\right|\\
    &\leq C\int_M|\overline{\ub}|
    \int_{-H}^0|\nabla\ub'||\ub'|^5\ud z\ud x\ud y.
\end{align*}
By the Cauchy--Schwarz inequality and the Holder's inequality, we find
\begin{multline}
   \left|\int_M\overline{\ub}
    \int_{-H}^0|\nabla\ub'||\ub'|^5\ud z\ud x\ud y\right|
    \leq \\
    \left(\int_M|\overline{\ub}|^4\ud M\right)^{\frac{1}{4}}
    \left(\int_\Omega|\nabla\ub'|^2|\ub'|^4\ud\Omega\right)^{\frac{1}{2}}
    \left(\int_M\left(\int_{-H}^0|\ub'|^6\ud z\right)^2\ud M\right)^{\frac{1}{4}}.
   \label{e4.19_a}
\end{multline}
By the Minkowski integral inequality \eqref{a6}, we have
\begin{equation}
   \left(\int_M\left(\int_{-H}^0|\ub'|^6\ud z\right)^2\ud M\right)^{\frac{1}{2}}
   \leq \int_{-H}^0\left(\int_M|\ub'|^{12}\ud M\right)^{\frac{1}{2}}\ud z.
   \label{e4.19_b}
\end{equation}
Applying the Ladyzhenskaya inequality \eqref{a2} to $\phi^3$ in $\mathbb{R}^2$, we obtain
\begin{equation}
   |\phi|_{L^{12}(M)}^{12} \leq C|\phi|_{L^{6}(M)}^{6}
   \left(\int_M|\phi|^4|\nabla\phi|^2\ud x\ud y\right) +
   |\phi|_{L^{6}(M)}^{12}.
   \label{e4.19a}
\end{equation}
Using \eqref{e4.19a}, we infer from \eqref{e4.19_b} that
\begin{align*}
   &\phantom{\leq}\left(\int_M\left(\int_{-H}^0|\ub'|^6
   \ud z\right)^2\ud M\right)^{\frac{1}{2}}\\
   &\leq \int_{-H}^0\left(C\int_M|\ub'|^{6}\ud M
   \int_M|\ub'|^4|\nabla\ub'|^2\ud M +
   \left(\int_M|\ub'|^6\ud M\right)^2 \right)^{\frac{1}{2}}\ud z\\
   &\leq C\int_{-H}^0\left(\int_M|\ub'|^{6}\ud M\right)^{\frac{1}{2}}
   \left(\int_M|\ub'|^4|\nabla\ub'|^2\ud M\right)^{\frac{1}{2}} +
   \int_\Omega|\ub'|^6\ud \Omega\\
   &\leq C\left(\int_\Omega|\ub'|^{6}\ud \Omega\right)^{\frac{1}{2}}
   \left(\int_\Omega|\ub'|^4|\nabla\ub'|^2\ud \Omega\right)^{\frac{1}{2}} +
   \int_\Omega|\ub'|^6\ud \Omega.
\end{align*}
Therefore we have
\begin{multline}
   \left(\int_M\left(\int_{-H}^0|\ub'|^6
      \ud z\right)^2\ud M\right)^{\frac{1}{2}} \leq \\
     C|\ub'|^{3}_{L^6}
   \left(\int_\Omega|\ub'|^4|\nabla\ub'|^2\ud \Omega\right)^{\frac{1}{2}} +
   |\ub'|^6_{L^6}.
   \label{e4.19b}
\end{multline}
By the Ladyzhenskaya inequality \eqref{a2} for functions in $\mathbb{R}^2$,
we find that
\begin{equation}
   \left(\int_M|\overline{\ub}|^4\ud M\right)^{\frac{1}{4}} \leq
   C|\overline{\ub}|_{L^2(M)}^{\frac{1}{2}}|\overline{\ub}|_{H^1(M)}^{\frac{1}{2}},
   \label{e4.19c}
\end{equation}
Using \eqref{e4.19b} and \eqref{e4.19c} we infer from \eqref{e4.19_a} that
\begin{multline}
   \int_M\left(\overline{\ub}
    \int_{-H}^0|\nabla\ub'||\ub'|^5\ud z\right)\ud x\ud y\\
    \leq C|\overline{\ub}|_{L^2(M)}^{\frac{1}{2}}
    |\overline{\ub}|_{H^1(M)}^{\frac{1}{2}}\left(
    |\ub'|_{L^6}^{\frac{3}{2}}\left(\int_\Omega|\ub'|^4
    |\nabla\ub'|^2\ud \Omega\right)^{\frac{3}{4}} +
    |\ub'|_{L^6}^6\right).
    \label{e4.19d}
\end{multline}
By Young's inequality, we have
\begin{multline}
   |\overline{\ub}|_{L^2(M)}^{\frac{1}{2}}
    |\overline{\ub}|_{H^1(M)}^{\frac{1}{2}}
    |\ub'|_{L^6}^{\frac{3}{2}}\left(\int_\Omega|\ub'|^4
    |\nabla\ub'|^2\ud \Omega\right)^{\frac{3}{4}}\\
   \leq C|\overline{\ub}|_{L^2(M)}^{2}
    |\overline{\ub}|_{H^1(M)}^{2}|\ub'|_{L^6}^6 +
    \dfrac{1}{4}\mu\int_\Omega|\ub'|^4
    |\nabla\ub'|^2\ud \Omega,
    \label{e4.19e}
\end{multline}
For the other integral that involves nonlinear terms,
\begin{align*}
   &\phantom{=} \left| \int_\Omega \overline{(\ub'\cdot\nabla)\ub' +
   (\nabla\cdot\ub')\ub'}\cdot|\ub'|^4\ub'\ud\Omega \right| \\
   &= \left|\int_\Omega\overline{u'_iu'_j}\p_i(|\ub'|^4u'_j)\ud\Omega\right|\\
    &\leq \left|\int_M\overline{u_i'u'_j}
     \left(\int_{-H}^0\p_j(|\ub'|^4u'_j)\ud z\right)\ud x\ud y\right|\\
    &\leq \left|\dfrac{1}{H}\int_M\left(\int_{-H}^0{u_i'u'_j}\ud z
     \int_{-H}^0\p_j(|\ub'|^4u'_j)\ud z\right)\ud x\ud y\right|\\
    &\leq C\int_M\left(\int_{-H}^0|\ub'|^2\ud z
     \int_{-H}^0|\ub'|^4|\nabla\ub'|^2\ud z\right)\ud x\ud y.
\end{align*}
By similar use of the Minkowski inequality \eqref{a6} and various interpolation inequalities,
we find that
\begin{multline}
   \int_M\left(\int_{-H}^0|\ub'|^2\ud z
     \int_{-H}^0|\ub'|^4|\nabla\ub'|^2\ud z\right)\ud x\ud y \leq \\
     C|\ub'|_{L^6(\Omega)}^3\left(\int_\Omega|\ub'|^4
     |\nabla\ub'|^2\ud \Omega\right)^{\frac{1}{2}}\left(
     |\ub'|_{L^2(\Omega)} + |\nabla\ub'|_{L^2(\Omega)}\right).
     \label{e4.19f}
\end{multline}
By Young's inequality, we have
\begin{multline}
     \left(|\ub'|_{L^2(\Omega)} + |\nabla\ub'|_{L^2(\Omega)}\right)
    |\ub'|_{L^6(\Omega)}^3\left(\int_\Omega|\ub'|^4
     |\nabla\ub'|^2\ud \Omega\right)^{\frac{1}{2}} \leq \\
   \left(|\ub'|^2_{L^2(\Omega)} + |\nabla\ub'|^2_{L^2(\Omega)}\right)
    |\ub'|_{L^6(\Omega)}^6 +
     \dfrac{1}{4}\mu\int_\Omega|\ub'|^4
    |\nabla\ub'|^2\ud \Omega.
    \label{e4.19g}
\end{multline}
Putting \eqref{e4.19}--\eqref{e4.19g} together, we find
\begin{multline}
   \dfrac{\ud}{\ud t}|\ub'|_{L^6}^6 +
   \mu\int_\Omega\left(|\nabla\ub'|^2|\ub'|^4 +
   \left|\nabla|\ub'|^2\right|^2|\ub'|^2\right)\ud\Omega +{ }\\
   \nu\int_\Omega\left(\left|\dfrac{\p \ub'}{\p z}\right|^2|\ub'|^4 +
   \left|\dfrac{\p}{\p z}|\ub'|^2\right|^2|\ub'|^2 \right)\ud\Omega 
   \leq
   \left(C|\mathbf{F}'|^2_{L^2} + C_2|\ub'|^2_{L^2}\right)
   |\ub'|_{L^6}^4 + { }\\
     C\left( |\overline{\ub}|_{L^2(M)}^{\frac{1}{2}}
    |\overline{\ub}|_{H^1(M)}^{\frac{1}{2}} +
    |\overline{\ub}|_{L^2(M)}^2|\overline{\ub}|_{H^1(M)}^2 +
    |\ub'|_{L^2(\Omega)}^2 + C|\nabla\ub'|_{L^2(\Omega)}^2 \right)
    |\ub'|^6_{L^6}.
    \label{e4.20}
\end{multline}
Ignoring the other positive terms on the left hand side of \eqref{e4.20},
and dividing both sides by $|\ub'|_{L^6}^4$,
we have
\begin{multline}
   \dfrac{\ud}{\ud t}|\ub'|_{L^6}^2 \leq
   C|\mathbf{F}'|^2_{L^2} + C_2|\ub'|_{L^2}^2 + \\
     C\left(|{\ub}|_{L^2(\Omega)}^2|{\ub}|_{H^1(\Omega)}^2 +
    |\ub|_{L^2(\Omega)}^2 + |\nabla\ub|_{L^2(\Omega)}^2 \right)
    |\ub'|^2_{L^6}.
    \label{e4.21}
\end{multline}
Applying the Gronwall inequality to \eqref{e4.21}, and using the $L^2$ estimate
result \eqref{e4.16}, we obtain
\begin{equation}
   |\ub'(\cdot,t)|_{L^6}^2 \leq J_6,
   \label{e4.22}
\end{equation}
with
\begin{displaymath}
   J_6(t) =
   e^{C\left(K_1^2(t) + K_1(t)t + K_1(t)\right) }\left(|\ub_0|_{L^6(\Omega)}^2 +
   C_2J_1(t) + C\int_0^t|\mathbf{F}|_{L^2(\Omega)}^2\ud s\right).
\end{displaymath}
Integrating \eqref{e4.20} over $[0,t]$, and using the estimate
\eqref{e4.22}, we obtain
\begin{multline}
   \mu\int_0^t\int_\Omega\left(|\nabla\ub'|^2|\ub'|^4 +
   \left|\nabla|\ub'|^2\right|^2|\ub'|^2\right)\ud\Omega\ud s +{ }\\
   \nu\int_0^t\int_\Omega\left(\left|\dfrac{\p \ub'}{\p z}\right|^2|\ub'|^4 +
   \left|\dfrac{\p}{\p z}|\ub'|^2\right|^2|\ub'|^2 \right)\ud\Omega\ud s
      \leq \tilde J_6(t),
    \label{e4.23}
\end{multline}
with
\begin{multline*}
   \tilde J_6(t) = |\ub_0|_{L^6}^6 +
   J^2_6(t)\left(C\int_0^t|F(\cdot,s)|_{L^2}^2\ud s + C_2K_1t\right) +{ }\\
   CJ^3_6(t)\left(K_1^2(t) + K_1(t)t + K_1(t) \right).
\end{multline*}

\vspace{5mm}
\noindent{\it Estimate $|\nabla\overline\ub|_{L^2(M)}$}\\
We multiply \eqref{e4.4a} by $-\Delta\overline{\ub}$ and
integrate by parts over $M$
to obtain
\begin{multline}
   \dfrac{1}{2}\dfrac{\ud}{\ud t}|\nabla\overline{\ub}|_{L^2(M)}^2 +
   \mu|\Delta\overline{\ub}|_{L^2(M)}^2 +
   \int_Mf\bk\times\overline{\ub}\cdot(-\Delta\overline{\ub})\ud M = \\
   \int_M\overline{\mathbf{F}}\cdot\Delta\overline{\ub}\ud M -
   \dfrac{1}{\rho_0}\int_M\nabla p_0\cdot\Delta\overline{\ub}\ud M -
   \int_M(\overline{\ub}\cdot\nabla)\overline{\ub}\cdot\Delta\overline{\ub}\ud M\\
   - \int_M\overline{(\ub'\cdot\nabla)\ub' + (\nabla\cdot\ub')\ub'}\cdot
   \Delta\overline{\ub}\ud M +
   \mu\int_\Omega\Delta P_{M,N} \overline{\ub}\cdot
   \Delta\overline{\ub}\ud\Omega.
   \label{e4.24}
\end{multline}
We note that, by integration by parts,
\begin{align*}
    \left|\int_M(f\bk\times \overline{\ub})\cdot\Delta\overline{\ub}\ud\Omega\right|
    &\leq |f|_\infty |\overline{\ub}|_{L^2}|\Delta\overline{\ub}|_{L^2}\\
    &\leq C|f|_\infty^2|\overline{\ub}|_{L^2}^2 + \dfrac{\mu}{4}|\Delta\overline{\ub}|_{L^2}^2.
\end{align*}
And thanks to \eqref{e4.4a6},
\begin{equation*}
   \int_M\nabla p_0\cdot\Delta\overline{\ub}\ud M = 0.
\end{equation*}
Following similar steps in the handling of the 2D Navier Stokes equations, we obtain
\begin{equation*}
   \int_M(\overline{\ub}\cdot\nabla)\overline{\ub}\cdot\Delta\overline{\ub}\ud M
   \leq C|\overline{\ub}|_{L^2(M)}^{\frac{1}{2}}
   |\nabla\overline{\ub}|_{L^2(M)}|\Delta\overline{\ub}|_{L^2(M)}^{\frac{3}{2}}.
\end{equation*}
Applying the Cauchy--Schwarz and Holder inequalities one have
\begin{multline*}
   \int_M\overline{(\ub'\cdot\nabla)\ub' + (\nabla\cdot\ub')\ub'}\cdot
   \Delta\overline{\ub}\ud M \leq \\
   C|\nabla\ub'|_{L^2(\Omega)}^{\frac{1}{2}}\left(
   \int_\Omega|\ub'|^4|\nabla\ub'|^2\ud\Omega\right)^{\frac{1}{4}}
   |\Delta\overline{\ub}|_{L^2(M)}.
\end{multline*}
For this last integral in \eqref{e4.24}, we proceed by the Cauchy--Schwarz inequality
and \eqref{e4.14} with $k=2$,
\begin{align*}
   \mu\int_\Omega\Delta P_{M,N} \overline{\ub}\cdot
   \Delta\overline{\ub}\ud\Omega &\leq
   \mu|\Delta P_{M,N}\overline{\ub}|_{L^2}
   |\Delta\overline{\ub}|_{L^2}^2\\
   &\leq
   C_2\mu|\overline{\ub}|_{L^2}|\Delta\overline\ub|_{L^2}^2\\
   &\leq
   C_2\mu|\overline{\ub}|_{L^2}^2 +
   \dfrac{\mu}{4}|\Delta\overline\ub|_{L^2}^2.\\
\end{align*}
Then we derive from \eqref{e4.24} that
\begin{multline}
   \dfrac{\ud}{\ud t}|\nabla\overline{\ub}|_{L^2(M)}^2 +
   \mu|\Delta\overline{\ub}|_{L^2(M)}^2 \leq
   \dfrac{2}{\mu}|\overline{\mathbf{F}}|_{L^2(M)}^2 +
   (C_2\mu+C|f|_\infty^2)|\ub|_{L^2}^2 + \\
   C|\overline{\ub}|_{L^2(M)}^2|\nabla\overline{\ub}|_{L^2(M)}^4 +
   C|\nabla\ub'|_{L^2(\Omega)}^2 + C\int_\Omega|\ub'|^4|\nabla\ub'|^2\ud\Omega.
    \label{e4.25}
\end{multline}
Applying the Gronwall inequality to \eqref{e4.25}, and using the previous estimates
\eqref{e4.16} and \eqref{e4.23}, we obtain
\begin{equation}
   |\nabla\overline{\ub}(\cdot,t)|_{L^2(M)}^2 +
   \dfrac{\mu}{2}\int_0^t|\Delta\overline{\ub}|_{L^2(M)}^2\ud s \leq
   J_2(t),
   \label{e4.26}
\end{equation}
with
\begin{multline*}
   J_2(t) = e^{CJ_1^2(t)}\Big(|\ub_0|_{H^1}^2 +
   (C_2\mu + C|f|_\infty^2) J_1(t)t + { } \\
   C\int_0^t|\mathbf{F}|_{L^2(M)}^2\ud s
   + CJ_1(t) + C\tilde J_6(t) \Big).
\end{multline*}

\vspace{5mm}
\noindent{\it To estimate $|\p_z\ub|_{L^2(\Omega)}$}\\
We multiply \eqref{e4.4a0} by $\p^2\ub/\p z^2$ and integrate by parts over $\Omega$,
and  we have
\begin{multline}
   \dfrac{1}{2}\dfrac{\ud}{\ud t}|\p_z\ub|_{L^2(\Omega)}^2 +
   \mu|\nabla\p_z\ub|_{L^2(\Omega)}^2 + \nu\left|\dfrac{\p^2\ub}{\p z^2}\right|_{L^2(\Omega)}^2
    = \left(F, \dfrac{\p^2\ub}{\p z^2}\right) -\\
   \int_\Omega\left( (\ub\cdot\nabla)\ub + \left(\nabla\cdot\int_z^0\ub\ud z\right)
   \dfrac{\p\ub}{\p z}\right)\cdot\dfrac{\p^2\ub}{\p z^2}\ud\Omega +{ }\\
   \mu|\nabla P_{M,N} \p_z\ub|_{L^2(\Omega)}^2 +
   \nu\left|\dfrac{\p^2}{\p z^2} P_{M,N} \ub\right|_{L^2(\Omega)}^2.
   \label{e4.27}
\end{multline}
By Holder's inequality,
\begin{equation*}
   \left(F, \dfrac{\p^2\ub}{\p z^2}\right) \leq
   C|F|_{L^2}^2 + \dfrac{\nu}{4}\left|\dfrac{\p^2\ub}{\p z^2}\right|_{L^2}^2.
\end{equation*}
With regard to the integral in \eqref{e4.27} that involves the nonlinear convection
terms, we find
\begin{align*}
   \phantom{=} & - \int_\Omega\left( (\ub\cdot\nabla)\ub +
   \left(\nabla\cdot\int_z^0\ub\ud z\right)
   \dfrac{\p\ub}{\p z}\right)\cdot\dfrac{\p^2\ub}{\p z^2}\ud\Omega\\
   = &\int_\Omega\dfrac{\p}{\p z}\left( (\ub\cdot\nabla)\ub +
   \left(\nabla\cdot\int_z^0\ub\ud z\right)
   \dfrac{\p\ub}{\p z}\right)\cdot\dfrac{\p\ub}{\p z}\ud\Omega\\
   = &\int_\Omega\left( (\p_z\ub\cdot\nabla)\ub + (\ub\cdot\nabla)\p_z\ub -
   (\nabla\cdot\ub)\dfrac{\p\ub}{\p z} + \nabla\cdot\int_z^0\ub\ud\xi
   \dfrac{\p^2\ub}{\p z^2}\right)\cdot\dfrac{\p\ub}{\p z}\ud\Omega\\
   = &\int_\Omega\left( (\p_z\ub\cdot\nabla)\ub  -
   (\nabla\cdot\ub)\dfrac{\p\ub}{\p z} \right)
   \cdot\dfrac{\p\ub}{\p z}\ud\Omega\\
   = &\int_\Omega -(\nabla\cdot\p_z\ub)\ub\cdot\p_z\ub -
      (\p_z\ub\cdot\nabla)\p_z\ub\cdot\ub +
      2(\ub\cdot\nabla)\p_z\ub\cdot\p_z\ub\ud\Omega.
\end{align*}
Hence
\begin{align*}
   \phantom{\leq}&\left| \int_\Omega\left( (\ub\cdot\nabla)\ub +
   \left(\nabla\cdot\int_z^0\ub\ud z\right)
   \dfrac{\p\ub}{\p z}\right)\cdot\dfrac{\p^2\ub}{\p z^2}\ud\Omega\right|\\
   &\leq C\int_\Omega |\ub||\nabla\p_z\ub||\p_z\ub|\ud\Omega \\
   &\leq C|\ub|_{L^6(\Omega)}|\p_z\ub|_{L^3(\Omega)}|\nabla\p_z\ub|_{L^2(\Omega)}\\
   &\leq C|\ub|_{L^6(\Omega)}|\p_z\ub|_{L^2(\Omega)}^{\frac{1}{2}}
   |\p_z\ub|_{H^1(\Omega)}^\frac{1}{2}|\nabla\p_z\ub|_{L^2(\Omega)}\\
   &\leq C|\ub|_{L^6(\Omega)}|\p_z\ub|_{L^2(\Omega)}^{\frac{1}{2}}
   |\ub_{zz}|_{L^2}^\frac{1}{2}|\nabla\p_z\ub|_{L^2(\Omega)} +
   C|\ub|_{L^6(\Omega)}|\p_z\ub|_{L^2(\Omega)}^{\frac{1}{2}}
   |\nabla\p_z\ub|_{L^2(\Omega)}^\frac{3}{2}\\
   &\leq C|\ub|_{L^6}^4|\p_z\ub|_{L^2}^2 + \dfrac{\nu}{4}|\ub_{zz}|_{L^2}^2 +
   \dfrac{\mu}{2}|\nabla\p_z\ub|_{L^2}^2.
\end{align*}
The last two integrals on the right--hand side of \eqref{e4.27}
can be handled by \eqref{e4.14} with $k=2$. Thus we derive from
\eqref{e4.27} that
\begin{multline}
   \dfrac{\ud}{\ud t}|\p_z\ub|_{L^2(\Omega)}^2 +
   \mu|\nabla\p_z\ub|_{L^2(\Omega)}^2 + \nu\left|\dfrac{\p^2\ub}{\p z^2}\right|_{L^2(\Omega)}^2
   \leq C|\mathbf{F}|_{L^2}^2 + { }\\
   C\left( J_2^2(t) + J_6^{2}(t)\right) |\p_z\ub|_{L^2}^2
   + C_2|\ub|_{L^2}^2.
   \label{e4.28}
\end{multline}
An application of the Gronwall inequality to \eqref{e4.28}
readily yields
\begin{multline}
   |\p_z\ub(\cdot,t)|_{L^2(\Omega)}^2 +
   \mu\int_0^t|\nabla\p_z\ub|_{L^2(\Omega)}^2\ud s +
   \nu\int_0^t|\dfrac{\p^2\ub}{\p z^2}|_{L^2(\Omega)}^2\ud s
    \leq J_z(t),
   \label{e4.29}
\end{multline}
with
\begin{displaymath}
    J_z(t) = e^{C\left( J_2^2(t) + J_6^{2}(t)\right)t}
   \left(|\ub_0|_{H^1}^2 + C\int_0^t|\mathbf{F}|_{L^2}^2\ud
    + C_2J_1 t\right).
\end{displaymath}

\vspace{5mm}
\noindent{\it To estimate $|\nabla\ub|_{L^2}$}\\
We multiply \eqref{e4.4a0} by $-\Delta\ub$ and integrate by parts over $\Omega$,
\begin{multline}
   \dfrac{1}{2}\dfrac{\ud}{\ud t}|\nabla\ub|_{L^2(\Omega)}^2 +
   \int_\Omega\left( (\ub\cdot\nabla)\ub + \left(\nabla\cdot\int_z^0\ub\ud\xi\right)
   \dfrac{\p\ub}{\p z}\right)\cdot(-\Delta\ub)\ud\Omega +{ } \\
   f\int_{\Omega}k\times\ub\cdot(-\Delta\ub)\ud\Omega +
   \dfrac{1}{\rho_0}\int_\Omega\nabla p_0\cdot(-\Delta\ub)\ud\Omega +
   \mu|\Delta\ub|_{L^2(\Omega)}^2 +
   \nu|\nabla\p_z\ub|_{L^2(\Omega)}^2 \\
     = \left(F, -\Delta\ub\right) +
     \mu|\Delta P_{M,N} \ub|_{L^2(\Omega)}^2 +
     \nu|\nabla P_{M,N} \p_z\ub|_{L^2(\Omega)}^2.
   \label{e4.30}
\end{multline}
We note that
\begin{displaymath}
    \int_M(f\bk\times {\ub})\cdot\Delta{\ub}\ud\Omega = 0,
\end{displaymath}
and
\begin{align*}
   \int_\Omega\nabla p_0\cdot(-\Delta\ub)\ud\Omega &=
   -\int_M\nabla p_0\cdot\int_{-H}^0(\Delta\ub)\ud z\ud M\\
   &= -{H}\int_M\nabla p_0\cdot\Delta\overline{\ub}\ud M\\
   & = 0.
\end{align*}
It is easy to see that
\begin{multline*}
   \int_\Omega\left( (\ub\cdot\nabla)\ub + \left(\nabla\cdot\int_z^0\ub\ud\xi\right)
   \dfrac{\p\ub}{\p z}\right)\cdot(\Delta\ub)\ud\Omega \leq \\
   C\int_\Omega\left(|\ub||\nabla\ub| + \int_{-H}^0|\nabla\ub|\ud z|\p_z\ub|\right)
   |\Delta\ub|\ud\Omega.
\end{multline*}
By the Young's inequality and the Sobolev interpolation inequality \eqref{a4} for $L^3$
functions, we have
\begin{align*}
   \int_\Omega |\ub||\nabla\ub||\Delta\ub|\ud\Omega
   &\leq C|\ub|_{L^6}|\nabla\ub|_{L^3}|\Delta\ub|_{L^2}\\
   &\leq C|\ub|_{L^6}|\nabla\ub|_{L^2}^{\frac{1}{2}}|\nabla\p_z\ub|_{L^2}^\frac{1}{2}
   |\Delta\ub|_{L^2} +
   C|\ub|_{L^6}|\nabla\ub|_{L^2}^{\frac{1}{2}}|\Delta\ub|_{L^2}^\frac{3}{2}.
\end{align*}
We appeal to the following inequality for functions in
$\mathbb{R}^3$ (the scaling is 2D due to the vertical integration.
For a proof, see \cite{CT03}),
\begin{displaymath}
   \int_\Omega\left(\int_{-H}^0|\nabla\ub|\ud z\right)|f||g|\ud\Omega \leq
   C|f|_{L^2}|\ub|_{H^1}^{\frac{1}{2}}|\ub|_{H^2}^{\frac{1}{2}}
   |g|_{L^2}^{\frac{1}{2}}|g|_{H^1}^{\frac{1}{2}},
\end{displaymath}
to obtain (by setting $f=\Delta\ub$, $g=\p_z\ub$)
\begin{displaymath}
   \int_\Omega\left( \int_{-H}^0|\nabla\ub|\ud z\right)|\p_z\ub|
   |\Delta\ub|\ud\Omega \leq
   C|\nabla\ub|_{L^2}^{\frac{1}{2}}|\p_z\ub|_{L^2}^{\frac{1}{2}}
   |\nabla\p_z\ub|_{L^2}^{\frac{1}{2}}|\Delta\ub|_{L^2}^{\frac{3}{2}}.
\end{displaymath}
The last two integrals in the equation above are handled by
\eqref{e4.14} with $k=2$.
After these intermediate steps we
derive from \eqref{e4.30} that
\begin{multline}
   \dfrac{\ud}{\ud t}|\nabla\ub|_{L^2}^2 + {\mu}|\Delta\ub|_{L^2}^2
   + \nu|\nabla\p_z\ub|_{L^2}^2  \leq
   \dfrac{4}{\mu}|\mathbf{F}|_{L^2}^2 + { }\\
   C\left(|\ub|_{L^6}^4 + |\p_z\ub|_{L^2}^2|\nabla\p_z\ub|_{L^2}^2\right)
   |\nabla\ub|_{L^2}^2 + C_2|\ub|_{L^2}^2.
   \label{e4.31}
\end{multline}
We appeal to the Gronwall inequality, and, with use of \eqref{e4.16},
\eqref{e4.22}, and \eqref{e4.29}, we obtain
\begin{equation}
   |\nabla\ub(\cdot,t)|_{L^2}^2 + {\mu}\int_0^t|\Delta\ub|_{L^2}^2\ud s
   + \nu\int_0^t|\nabla\p_z\ub|_{L^2}^2\ud s
   \leq
   J_V(t),
   \label{e4.32}
\end{equation}
with
\begin{displaymath}
   J_V(t) =
   e^{C\left(J_6^{2} + J_2^2(t)\right)t + J_z^2(t)}\left(
   |\ub_0|_{H^1}^2 + C_2J_1 t + \dfrac{4}{\mu}\int_0^t|\mathbf{F}(\cdot,s)|_{L^2}^2\ud s
   \right).
\end{displaymath}

We now prove Theorem \ref{t0}.
\begin{proof}
    The short time existence and
    uniqueness of the strong solutions of \eqref{e4.1}--\eqref{e4.3e}
    can be established as for
    the viscous primitive equations (see \cite{GMR01, HTZ03}).
   Let $\ub$ be such a strong solution
   corresponding to the initial data $\ub_0$ with the
   maximal interval of existence $[0,T^*)$. If
   $T^* \geq T$, then there is nothing to prove here. Let
   us suppose that $T^*<T$, and in particular, $T^* <\infty$.
   Then it is clear that
   \begin{equation}
      \limsup_{t\rightarrow T^*-}||\ub||_{H^1} = \infty.
      \label{e4.32a}
   \end{equation}
   Otherwise the solution can be extended beyond $T^*$.
   However the estimates \eqref{e4.29} and \eqref{e4.32} indicate that
   $||\ub(\cdot,t)||_{H^1} < \infty$ for all $t<T$, which
   contradicts \eqref{e4.32a}. Hence the solution
   must exist for the whole period of $[0,T)$.

It remains to show the continuous dependence of the solution
on the data, of which the uniqueness of the solution is a consequence.
Let us assume that $\ub^1$ and $\ub^2$ be two solutions
   corresponding to the two sets of initial data
   $\ub^1_0$ and $\ub^2_0$, respectively. We let
   \begin{displaymath}
      \ub = \ub^1 - \ub^2.
   \end{displaymath}
   Then $\ub$ satisfies the following equation
   \begin{multline}
      \dfrac{\p\ub}{\p t} + (\ub^1\cdot\nabla)\ub +
      (\nabla\cdot\int_z^0\ub^1\ud\xi)\dfrac{\p\ub}{\p z} +
      (\ub\cdot\nabla)\ub^2 +
      (\nabla\cdot\int_z^0\ub\ud\xi)\dfrac{\p\ub^2}{\p z} \\
      +f\bk\times\ub  -
      \mu\Delta\ub
      - \nu\dfrac{\p^2\ub}{\p z^2} =
      -\mu\Delta P_{M,N} \ub
      - \nu\dfrac{\p^2}{\p z^2} P_{M,N} \ub.
      \label{e4.33}
   \end{multline}
   We multiply \eqref{e4.33} by $\ub$ and integrate by parts
   over $\Omega$ to obtain
   \begin{multline}
      \dfrac{1}{2}\dfrac{\ud}{\ud t}|\ub|_{L^2}^2 +
      \left( (\ub^1\cdot\nabla)\ub,\ub \right) + { }\\
      \left( (\nabla\cdot\int_z^0\ub^1\ud\xi)\dfrac{\p\ub}{\p z},\,\ub \right) +
      \left( (\ub\cdot\nabla)\ub^2,\ub \right) + 
      \left( (\nabla\cdot\int_z^0\ub\ud\xi)\dfrac{\p\ub^2}{\p z},\,\ub \right) + { }\\
      \mu|\nabla\ub|^2 + \nu\left|\dfrac{\p \ub}{\p z}\right|^2  =
      \mu|\nabla P_{M,N} \ub|^2 +
      \nu\left|\dfrac{\p }{\p z} P_{M,N} \ub\right|^2.
      \label{e4.34}
   \end{multline}
   We verify that
   \begin{equation}
      \left( (\ub^1\cdot\nabla)\ub,\ub \right) +
      \left( (\nabla\cdot\int_z^0\ub^1\ud\xi)\dfrac{\p\ub}{\p z},\,\ub \right)
      =0,
      \label{e4.34a}
   \end{equation}
   \begin{equation}
      \left|\left( (\ub\cdot\nabla)\ub^2,\ub \right)\right|
      \leq |\nabla\ub^2|_{L^2}|\ub|_{L^2}^{\frac{1}{2}}
      |\nabla\ub|_{L^2}^{\frac{3}{2}},
      \label{e4.34b}
   \end{equation}
   and that
   \begin{equation}
      \left( (\nabla\cdot\int_z^0\ub\ud\xi)\dfrac{\p\ub^2}{\p z},\,\ub \right)
      \leq C|\ub|_{L^2}^{\frac{1}{2}}|\nabla\ub|_{L^2}^{\frac{3}{2}}
      |\p_z\ub^2|_{L^2}^{\frac{1}{2}}|\nabla\p_z\ub^2|_{L^2}^{\frac{1}{2}}.
      \label{e4.34c}
   \end{equation}
   Again, the last two integrals on the right--hand side of \eqref{e4.34}
   are handled by \eqref{e4.14} with $k=1$.
   Thus we derive from \eqref{e4.34} that
   \begin{multline}
      \dfrac{\ud}{\ud t}|\ub|_{L^2}^2 +
      \mu|\nabla\ub|^2 + \nu\left|\dfrac{\p \ub}{\p z}\right|^2 
      \leq \\
      C\left(|\nabla\ub^2|_{L^2}^4 + |\p_z\ub^2|_{L^2}^2|\nabla\p_z\ub^2|_{L^2}^2
      \right)|\ub|_{L^2}^2 + C_1|\ub|_{L^2}^2.
      \label{e4.35}
   \end{multline}
   Thanks to the a priori estimates \eqref{e4.29} and \eqref{e4.32}
   and the Gronwall inequality, we have
   \begin{equation}
      |\ub(\cdot,t)|_{L^2}^2 \leq
      e^{C(K_V^2(t)t + K_z^2(t)) + C_1}|\ub(\cdot,0)|_{L^2}^2.
      \label{e4.36}
   \end{equation}
   This shows that the solution depends on the initial data continuously.
   When $\ub(\cdot,0) = \ub^1(\cdot,0)-\ub^2(\cdot,0) = 0$,
   \begin{equation}
      \ub(\cdot,t) = 0,\qquad\textrm{for all } t>0.
      \label{e4.37}
   \end{equation}
   This shows the uniqueness of the solution.
\end{proof}

\section{A spectral-viscosity model for geophysical turbulence}\label{s2}
In this section we study a model that has applications in the simulation
of geophysical turbulent flows.
\subsection{The model}\label{s2.1}
The 3D primitive equations with linear spectral eddy viscosity read
\begin{align}
   \dfrac{\p\ub}{\p t} + (\ub\cdot\nabla)\ub + w\dfrac{\p\ub}{\p z}
   +f\mathbf{k}\times\ub + \dfrac{1}{\rho_0}\nabla p - \mu\Delta\ub
   - \nu\dfrac{\p^2\ub}{\p z^2} - L\ub = \mathbf{F},\label{e1}\\
   \dfrac{\p p}{\p z} = -\rho_0g,\label{e2}\\
   \nabla\cdot\ub + \dfrac{\p w}{\p z} = 0,\label{e3}
\end{align}
with
\begin{equation*}
   L\ub = \mu_\delta\Delta (I-P_{M,N})\ub +
   \nu_\delta\dfrac{\p^2}{\p z^2}(I-P_{M,N})\ub.
\end{equation*}
The other notations being the same as those in Section \ref{s3},
the newly introduced ones $\mu_\delta$ and $\nu_\delta$ are the eddy
closure parameters in the horizontal and vertical directions,
respectively. The subscript $\delta$ indicates that the parameters
depend on the grid resolution.

As for the model with partial viscosity in Section \ref{s3},
we consider a rectangular domain $\Omega = M\times(-H,0)$, with
$M=(0,L_x)\times(0,L_y)$, and we also consider periodic boundary conditions
in both the $x$ and $y$ directions, and free--slip, non--penetration
boundary conditions in the vertical directions. More precisely,
\begin{align}
   \ub(x,y,z,t) = \ub(x+L_x,y,z,t),\label{e4}\\
   \ub(x,y,z,t) = \ub(x,y+L_y,z,t),\label{e5}\\
   \dfrac{\p\ub}{\p n}|_{z=0} =
   \dfrac{\p\ub}{\p n}|_{z=-H} = 0,\label{e6}\\
   w|_{z=0} = w|_{z=-H} = 0,\label{e7}\\
   p|_{z=0} = p_0.\label{e8}
\end{align}
As a consequence of the boundary conditions taken here, the spectral
low--pass filter $P_{M,N}$ can be, and is defined as in \eqref{e4.3f}.

For the model \eqref{e1}--\eqref{e8} we obtain two results.
The first is the global well--posedness of the model, and the
second is concerned with the convergence of the solutions of
\eqref{e1}--\eqref{e8} to those of the viscous primitive equations
as $\mu_\delta$ and $\nu_\delta$ tend to zero. The first result
shall come as no surprise, and therefore its proof will only be
briefly sketched. The second result will be discussed in more details.

As usual, $p$ and $w$ can be expressed in terms of $\ub$. Specifically,
integrating \eqref{e2} from $z$ to 0, and using the boundary condition
\eqref{e8} we obtain
\begin{equation}
   p(x,y,z,t) = p_0(x,y,t) - \rho_0gz,\label{e10}
\end{equation}
Integrating \eqref{e3} from $z$ to 0, and using the boundary condition
\eqref{e7}, we obtain
\begin{equation}
   w(x,y,z,t) = \int_z^0\nabla\cdot\ub\ud\xi = \nabla\cdot\int_z^0\ub\ud\xi.
   \label{e11}
\end{equation}
Setting $z=-H$ in \eqref{e11}, and using \eqref{e7} again, we find
\begin{equation}
   \nabla\cdot\int_{-H}^0\ub\ud z = 0.
   \label{e12}
\end{equation}
We substitute \eqref{e10} and \eqref{e11} into \eqref{e1}, and
obtain a single closed equation for $\ub$:
\begin{multline}
   \dfrac{\p\ub}{\p t} + (\ub\cdot\nabla)\ub +
   (\nabla\cdot\int_z^0\ub\ud\xi)\dfrac{\p\ub}{\p z}
   +f\bk\times\ub + \dfrac{1}{\rho_0}\nabla p_0 - \\
   \mu\Delta\ub
   - \nu\dfrac{\p^2\ub}{\p z^2} - L\ub = \mathbf{F}.
   \label{e12a}
\end{multline}
The prognostic variable $\ub$ satisfies the equation \eqref{e1},
the constraint \eqref{e12}, and
the boundary conditions \eqref{e4}--\eqref{e6}. To complete the
system, we also require $\ub$ to satisfy the following initial
condition:
\begin{equation}
   \ub(x,y,z,0) = \ub_0(x,y,z).
   \label{e13}
\end{equation}
The diagnostic variables $w$ and $p$ are given by \eqref{e11} and
\eqref{e10} respectively.

We now specify the barotropic and baroclinic modes of \eqref{e12a}.
We let
\begin{equation}
   \overline{\ub}(x,y,t) = \dfrac{1}{H}\int_{-H}^0\ub(x,y,z,t)\ud z,
   \label{e14}
\end{equation}
which denotes the barotropic mode of the primitive variables. We
also denote by
\begin{equation}
   \ub'(x,y,z,t) = \ub(x,y,z,t)-\overline{\ub}(x,y,t)
   \label{e15}
\end{equation}
the baroclinic mode of the primitive variables.

As in Section \ref{s3}, we derive the equations and conditions that the
barotropic mode $\overline{\ub}$ and the baroclinic mode $\ub'$ must
satisfy. They are summarized as follows.
The equation and boundary conditions for $\overline{\ub}$ are
\begin{gather}
   \dfrac{\p\overline{\ub}}{\p t} +
   (\overline{\ub}\cdot\nabla)\overline{\ub} +
   [\overline{(\ub'\cdot\nabla)\ub'} +
   \overline{(\nabla\cdot\ub')\ub'} ]
   +f\bk\times\overline{\ub} + \dfrac{1}{\rho_0}\nabla p_0 -
   \mu\Delta\overline{\ub} =
   \overline{\mathbf{F}},\label{e26}\\
   \overline{\ub}(x,y,t) = \overline{\ub}(x+L_x,y,t),\label{e27}\\
   \overline{\ub}(x,y,t) = \overline{\ub}(x,y+L_y,t).\label{e28}\\
   \nabla\cdot\overline{\ub} = 0. \label{e29}
\end{gather}
The equation, boundary conditions, and constraint for $\ub'$ are
\begin{align}
   &\dfrac{\p\ub'}{\p t} +
   (\ub'\cdot\nabla)\ub' +
   (\nabla\cdot\int_z^0\ub'\ud\xi)\dfrac{\p\ub'}{\p z}
   + f\bk\times\ub' - \mu\Delta\ub' - \nu\dfrac{\p^2\ub'}{\p z^2}
   - L\ub'\nonumber\\
   &\phantom{ \dfrac{\p\ub'}{\p t} +}+ \left[(\ub'\cdot\nabla)\overline{\ub} + (\overline{\ub}\cdot\nabla)\ub'
   - \left(\overline{(\ub'\cdot\nabla)\ub'}+
   \overline{(\nabla\cdot\ub')\ub'}\right)\right]
   = \mathbf{F}',\label{e30}\\
   &\ub'(x,y,t) = \ub'(x+L_x,y,t),\label{e30a}\\
   &\ub'(x,y,t) = \ub'(x,y+L_y,t),\label{e31}\\
   &\left.\dfrac{\p\ub'}{\p z}\right|_{z=0} = \left.\dfrac{\p \ub'}{\p z}\right|_{z=-H} = 0,
   \label{e32}\\
   &\int_{-H}^0 \ub'\ud z = 0.\label{e33}
\end{align}

\subsection{Existence and uniqueness of strong solutions}\label{s2.2}
The functional settings are the same as specified in Section \ref{s3}.
A priori estimates are the essential ingredients in the proof
of global well--posedness of \eqref{e1}--\eqref{e8}. We shall first
derive the a priori estimates for $\ub,\,\overline{\ub}$ and $\ub'$.

\vspace{5mm}
\noindent{\it $L^2$ estimates}\\
We take inner product of \eqref{e12a} with $\ub$,
integrate by parts over $\Omega$, using the boundary
conditions \eqref{e4}--\eqref{e8}, we obtain
\begin{multline}
   (\dfrac{\p\ub}{\p t},\ub) + \left( (\ub\cdot\nabla)\ub,\ub \right) +
   \left( (\nabla\cdot\int_z^0\ub\ud\xi)\dfrac{\p\ub}{\p z},\,\ub \right) + { }\\
   (f\bk\times\ub,\ub) + \dfrac{1}{\rho_0}(\nabla p_0,\ub) +
   \mu|\nabla\ub|^2 + \nu\left|\dfrac{\p \ub}{\p z}\right|^2 + { }\\
   \mu_\delta|\nabla (I-P_{M,N})\ub|^2 + 
   \nu_\delta\left|\dfrac{\p}{\p z}(I-P_{M,N})\ub\right|^2 =
   (\mathbf{F},\ub).
   \label{e34a}
\end{multline}
Following similar steps in Section \ref{s3}, we obtain that
\begin{multline}
   |\ub(\cdot,t)|^2 + \mu\int_0^t|\nabla\ub|^2\ud s +
   \nu\int_0^t\left|\dfrac{\p}{\p z}\ub\right|^2\ud s +
   \mu_\delta\int_0^t|\nabla (I-P_{M,N})\ub|^2\ud s +\\
   \nu_\delta\int_0^t\left|\dfrac{\p}{\p z}(I-P_{M,N})\ub\right|^2\ud s
   \leq K_1(t),
   \label{e37}
\end{multline}
where
\begin{equation}
   K_1(t) \equiv |u_0|^2 + C\int_0^t|\mathbf{F}|^2\ud s.
   \label{e38}
\end{equation}
We note that $K_1(t)$ is a non--decreasing positive function
in $t$. This notion will be useful later in the paper.

\vspace{5mm}
\noindent{\it $L^6$ estimate on $\ub'$}\\
We rewrite \eqref{e30} as
\begin{multline}
   \dfrac{\p\ub'}{\p t} +
   (\ub'\cdot\nabla)\ub' +
   \left(\nabla\cdot\int_z^0\ub'\ud\xi\right)\dfrac{\p\ub'}{\p z}
   + f\bk\times\ub' - (\mu+\mu_\delta)\Delta\ub' -
   (\nu+\nu_\delta)\dfrac{\p^2\ub'}{\p z^2}\\
   { }+\left[(\ub'\cdot\nabla)\overline{\ub} + (\overline{\ub}\cdot\nabla)\ub'
   - \left(\overline{(\ub'\cdot\nabla)\ub'}+
   \overline{(\nabla\cdot\ub')\ub'}\right)\right]
   = \mathbf{F}'- \\
   \mu_\delta\Delta P_{M,N} \ub -
   \nu_\delta\dfrac{\p^2}{\p z^2} P_{M,N} \ub.
   \label{e38a}
\end{multline}
Now \eqref{e38a} is in the same form as equation \eqref{e4.8},
and therefore the same techniques from Section \eqref{s3} can
be applied to yield the $L^6$ estimates on $\ub'$. We
omit the details and state the results as follows.
\begin{equation}
   |\ub'(\cdot,t)|_{L^6}^2 \leq K_6,
   \label{e68c}
\end{equation}
with
\begin{displaymath}
   K_6(t) =
   e^{CK_1^2(t) + CK_1(t)t + CK_1(t) + C_1(t)}\left(|\ub_0|_{L^6(\Omega)}^2 +
   C\int_0^t|\mathbf{F}|_{L^2(\Omega)}^2\ud s\right).
\end{displaymath}
We also have
\begin{multline}
   (\mu+\mu_\delta)\int_0^t\int_\Omega\left(|\nabla\ub'|^2|\ub'|^4 +
   \left|\nabla|\ub'|^2\right|^2|\ub'|^2\right)\ud\Omega\ud s +{ }\\
   (\nu+\nu_\delta)\int_0^t\int_\Omega\left(\left|\dfrac{\p \ub'}{\p z}\right|^2|\ub'|^4 +
   \left|\dfrac{\p}{\p z}|\ub'|^2\right|^2|\ub'|^2 \right)\ud\Omega\ud s
   \leq \tilde K_6(t),
    \label{e69}
\end{multline}
with
\begin{multline*}
   \tilde K_6(t) = |\ub_0|_{L^6}^6 +
   K^2_6(t)\int_0^t|F(\cdot,s)|_{L^2}^2\ud s + { }\\
   K^3_6(t)\left(CK_1^2(t) + CK_1(t)t + CK_1(t) + C_1(t)\right).
\end{multline*}

\vspace{5mm}
\noindent{\it Estimate $|\nabla\overline{\ub}|_{L^2(M)}$}\\
We multiply \eqref{e26} by $-\Delta\overline{\ub}$ and integrate by parts over $M$
to obtain
\begin{multline}
   \dfrac{1}{2}\dfrac{\ud}{\ud t}|\nabla\overline{\ub}|_{L^2(M)}^2 +
   \mu|\Delta\overline{\ub}|_{L^2(M)}^2 -
   \int_Mf\bk\times\overline{\ub}\cdot\Delta\ub\ud M = { }\\
   \int_M\overline{\mathbf{F}}\cdot\Delta\overline{\ub}\ud M - 
   \dfrac{1}{\rho_0}\int_M\nabla p_0\cdot\Delta\overline{\ub}\ud M -
   \int_M(\overline{\ub}\cdot\nabla)\overline{\ub}\cdot\Delta\overline{\ub}\ud M
   - \\
   \int_M\overline{(\ub'\cdot\nabla)\ub' + (\nabla\cdot\ub')\ub'}\cdot
   \Delta\overline{\ub}\ud M.
   \label{e70}
\end{multline}
We note that
\begin{equation*}
   \int_Mf\bk\times\overline{\ub}\cdot\Delta\ub\ud M = 0,
\end{equation*}
and
\begin{equation*}
   \int_M\nabla p_0\cdot\Delta\overline{\ub}\ud M = 0.
\end{equation*}
Following similar steps in the handling of the 2D Navier Stokes equations, we obtain
\begin{equation*}
   \int_M(\overline{\ub}\cdot\nabla)\overline{\ub}\cdot\Delta\overline{\ub}\ud M
   \leq C|\overline{\ub}|_{L^2(M)}^{\frac{1}{2}}
   |\nabla\overline{\ub}|_{L^2(M)}|\Delta\overline{\ub}|_{L^2(M)}^{\frac{3}{2}}.
\end{equation*}
Applying the Cauchy--Schwarz and Holder inequalities one have
\begin{multline*}
   \int_M\overline{(\ub'\cdot\nabla)\ub' + (\nabla\cdot\ub')\ub'}\cdot
   \Delta\overline{\ub}\ud M \leq \\
   C|\nabla\ub'|_{L^2(\Omega)}^{\frac{1}{2}}\left(
   \int_\Omega|\ub'|^4|\nabla\ub'|^2\ud\Omega\right)^{\frac{1}{4}}
   |\Delta\overline{\ub}|_{L^2(M)}.
\end{multline*}
Utilizing these estimates, and using Youngs inequality again, we derive from
\eqref{e70} that
\begin{multline}
   \dfrac{\ud}{\ud t}|\nabla\overline{\ub}|_{L^2(M)}^2 +
   \dfrac{\mu}{2}|\Delta\overline{\ub}|_{L^2(M)}^2 \leq
   \dfrac{2}{\mu}|\overline{\mathbf{F}}|_{L^2(M)}^2 + { }\\
   C|\overline{\ub}|_{L^2(M)}^2|\nabla\overline{\ub}|_{L^2(M)}^4 +
   C|\nabla\ub'|_{L^2(\Omega)}^2 + C\int_\Omega|\ub'|^4|\nabla\ub'|^2\ud\Omega.
    \label{e71}
\end{multline}
Applying the Gronwall inequality to \eqref{e71}, and using the previous estimates
\eqref{e37} and \eqref{e69}, we obtain
\begin{equation}
   |\nabla\overline{\ub}(\cdot,t)|_{L^2(M)}^2 +
   \dfrac{\mu}{2}\int_0^t|\Delta\overline{\ub}|_{L^2(M)}^2\ud s \leq
   K_2(t),
   \label{e72}
\end{equation}
with
\begin{equation*}
   K_2(t) = e^{CK_1^2(t)}\left(|\ub_0|_{H^1}^2 + C\int_0^t|\mathbf{F}|_{L^2(M)}^2\ud s
   + CK_1(t) + C\tilde K_6(t) \right).
\end{equation*}

\vspace{5mm}
\noindent{\it To estimate $|\p_z\ub|_{L^2(\Omega)}$}\\
We multiply \eqref{e12a} by $\p^2\ub/\p z^2$ and integrate by parts over $\Omega$,
and utilizing the boundary conditions \eqref{e4}--\eqref{e6}, we have
\begin{multline}
   \dfrac{1}{2}\dfrac{\ud}{\ud t}|\p_z\ub|_{L^2(\Omega)}^2 +
   \mu|\nabla\p_z\ub|_{L^2(\Omega)}^2 + \nu\left|\dfrac{\p^2\ub}{\p z^2}\right|_{L^2(\Omega)}^2
   + { }\\
   \mu_\delta|\nabla \p_z(I-P_{M,N})\ub|_{L^2(\Omega)}^2 +
   \nu_\delta\left|\dfrac{\p^2(I-P_{M,N})\ub}{\p z^2}\right|_{L^2(\Omega)}^2
   = \left(F, \dfrac{\p^2\ub}{\p z^2}\right) - { }\\
   \int_\Omega\left( (\ub\cdot\nabla)\ub + \left(\nabla\cdot\int_z^0\ub\ud z\right)
   \dfrac{\p\ub}{\p z}\right)\cdot\dfrac{\p^2\ub}{\p z^2}\ud\Omega.
   \label{e75}
\end{multline}
By Holder's inequality,
\begin{equation*}
   \left(F, \dfrac{\p^2\ub}{\p z^2}\right) \leq
   C|F|_{L^2}^2 + \dfrac{\nu}{4}\left|\dfrac{\p^2\ub}{\p z^2}\right|_{L^2}^2.
\end{equation*}
With regard to the last integral in \eqref{e75} we find
\begin{align*}
   \phantom{=} & - \int_\Omega\left( (\ub\cdot\nabla)\ub +
   \left(\nabla\cdot\int_z^0\ub\ud z\right)
   \dfrac{\p\ub}{\p z}\right)\cdot\dfrac{\p^2\ub}{\p z^2}\ud\Omega\\
   = &\int_\Omega\dfrac{\p}{\p z}\left( (\ub\cdot\nabla)\ub +
   \left(\nabla\cdot\int_z^0\ub\ud z\right)
   \dfrac{\p\ub}{\p z}\right)\cdot\dfrac{\p\ub}{\p z}\ud\Omega\\
   = &\int_\Omega\left( (\p_z\ub\cdot\nabla)\ub + (\ub\cdot\nabla)\p_z\ub -
   (\nabla\cdot\ub)\dfrac{\p\ub}{\p z} + \nabla\cdot\int_z^0\ub\ud\xi
   \dfrac{\p^2\ub}{\p z^2}\right)\cdot\dfrac{\p\ub}{\p z}\ud\Omega\\
   = &\int_\Omega\left( (\p_z\ub\cdot\nabla)\ub  -
   (\nabla\cdot\ub)\dfrac{\p\ub}{\p z} \right)
   \cdot\dfrac{\p\ub}{\p z}\ud\Omega\\
   = &\int_\Omega -(\nabla\cdot\p_z\ub)\ub\cdot\p_z\ub -
      (\p_z\ub\cdot\nabla)\p_z\ub\cdot\ub +
      2(\ub\cdot\nabla)\p_z\ub\cdot\p_z\ub\ud\Omega.
\end{align*}
Hence
\begin{align*}
   \phantom{\leq}&\left| \int_\Omega\left( (\ub\cdot\nabla)\ub +
   \left(\nabla\cdot\int_z^0\ub\ud z\right)
   \dfrac{\p\ub}{\p z}\right)\cdot\dfrac{\p^2\ub}{\p z^2}\ud\Omega\right|\\
   &\leq C\int_\Omega |\ub||\nabla\p_z\ub||\p_z\ub|\ud\Omega \\
   &\leq C|\ub|_{L^6(\Omega)}|\p_z\ub|_{L^3(\Omega)}|\nabla\p_z\ub|_{L^2(\Omega)}\\
   &\leq C|\ub|_{L^6(\Omega)}|\p_z\ub|_{L^2(\Omega)}^{\frac{1}{2}}
   |\p_z\ub|_{H^1(\Omega)}^\frac{1}{2}|\nabla\p_z\ub|_{L^2(\Omega)}\\
   &\leq C|\ub|_{L^6(\Omega)}|\p_z\ub|_{L^2(\Omega)}^{\frac{1}{2}}
   |\ub_{zz}|_{L^2}^\frac{1}{2}|\nabla\p_z\ub|_{L^2(\Omega)} +
   C|\ub|_{L^6(\Omega)}|\p_z\ub|_{L^2(\Omega)}^{\frac{1}{2}}
   |\nabla\p_z\ub|_{L^2(\Omega)}^\frac{3}{2}\\
   &\leq C|\ub|_{L^6}^4|\p_z\ub|_{L^2}^2 + \dfrac{\nu}{4}|\ub_{zz}|_{L^2}^2 +
   \dfrac{\mu}{2}|\nabla\p_z\ub|_{L^2}^2.
\end{align*}
After these intermediate steps, we derive from \eqref{e75} that
\begin{multline}
   \dfrac{\ud}{\ud t}|\p_z\ub|_{L^2(\Omega)}^2 +
   \mu|\nabla\p_z\ub|_{L^2(\Omega)}^2 + \nu\left|\dfrac{\p^2\ub}{\p z^2}\right|_{L^2(\Omega)}^2
   + 2\mu_\delta|\nabla \p_z(I-P_{M,N})\ub|_{L^2(\Omega)}^2 + { }\\
   2\nu_\delta\left|\dfrac{\p^2(I-P_{M,N})\ub}{\p z^2}\right|_{L^2(\Omega)}^2
   \leq C|F|_{L^2}^2 +
   C|\ub|_{L^6}^4 |\p_z\ub|_{L^2}^2.
   \label{e76}
\end{multline}
We note that
\begin{align*}
   |\ub|_{L^6} &= |\overline{\ub} + \ub'|_{L^6}\\
   &\leq |\overline{\ub}|_{L^6} + |\ub'|_{L^6} \\
   &\leq C|\nabla\overline{\ub}|_{L^2} + |\ub'|_{L^6}.
\end{align*}
Therefore, by \eqref{e68c} and \eqref{e72}, we have
\begin{equation}
   |\ub|_{L^6}^4 \leq C\left( K_2^2(t) + K_6^{2}(t)\right).
   \label{e77}
\end{equation}
Then \eqref{e76}, together with \eqref{e77}, gives
\begin{multline}
   \dfrac{\ud}{\ud t}|\p_z\ub|_{L^2(\Omega)}^2 +
   \mu|\nabla\p_z\ub|_{L^2(\Omega)}^2 + \nu\left|\dfrac{\p^2\ub}{\p z^2}\right|_{L^2(\Omega)}^2
   + 2\mu_\delta|\nabla \p_z(I-P_{M,N})\ub|_{L^2(\Omega)}^2 + { }\\
   2\nu_\delta\left|\dfrac{\p^2(I-P_{M,N})\ub}{\p z^2}\right|_{L^2(\Omega)}^2
   \leq C|\mathbf{F}|_{L^2}^2 +
   C\left( K_2^2(t) + K_6^{2}(t)\right) |\p_z\ub|_{L^2}^2.
   \label{e77a}
\end{multline}
An application of the Gronwall inequality to \eqref{e77a} readily gives
\begin{multline}
   |\p_z\ub(\cdot,t)|_{L^2(\Omega)}^2 +
   \mu\int_0^t|\nabla\p_z\ub|_{L^2(\Omega)}^2\ud s + { }\\
   \nu\int_0^t|\dfrac{\p^2\ub}{\p z^2}|_{L^2(\Omega)}^2\ud s 
   + 2\mu_\delta\int_0^t|\nabla \p_z(I-P_{M,N})\ub|_{L^2(\Omega)}^2\ud s + { }\\
   2\nu_\delta\int_0^t |\dfrac{\p^2(I-P_{M,N})\ub}{\p z^2}|_{L^2(\Omega)}^2\ud s
   \leq K_z(t),
   \label{e78}
\end{multline}
with
\begin{displaymath}
    K_z(t) = e^{C\left( K_2^2(t) + K_6^{2}(t)\right)t}
   \left(|\ub_0|_{H^1}^2 + C\int_0^t|\mathbf{F}|_{L^2}^2\ud s\right).
\end{displaymath}

\vspace{5mm}
\noindent{\it To estimate $|\nabla\ub|_{L^2}$}\\
We multiply \eqref{e12a} by $-\Delta\ub$ and integrate by parts over $\Omega$,
\begin{multline}
   \dfrac{1}{2}\dfrac{\ud}{\ud t}|\nabla\ub|_{L^2(\Omega)}^2 +
   \int_\Omega\left( (\ub\cdot\nabla)\ub + \left(\nabla\cdot\int_z^0\ub\ud\xi\right)
   \dfrac{\p\ub}{\p z}\right)\cdot(-\Delta\ub)\ud\Omega +{ } \\
   f\int_{\Omega}k\times\ub\cdot(-\Delta\ub)\ud\Omega +
   \dfrac{1}{\rho_0}\int_\Omega\nabla p_0\cdot(-\Delta\ub)\ud\Omega +
   \mu|\Delta\ub|_{L^2(\Omega)}^2 +{ } \\
   \nu|\nabla\p_z\ub|_{L^2(\Omega)}^2
   + \mu_\delta|\Delta (I-P_{M,N})\ub|_{L^2(\Omega)}^2 +
   \nu_\delta|\nabla \p_z(I-P_{M,N})\ub|_{L^2(\Omega)}^2\\
   = \left(F, -\Delta\ub\right)
   \label{e79}
\end{multline}
We note that
\begin{displaymath}
    \int_M(f\bk\times {\ub})\cdot\Delta{\ub}\ud\Omega = 0,
\end{displaymath}
and
\begin{align*}
   \int_\Omega\nabla p_0\cdot(-\Delta\ub)\ud\Omega &=
   -\int_M\nabla p_0\cdot\int_{-H}^0(\Delta\ub)\ud z\ud M \\
   &= -{H}\int_M\nabla p_0\cdot\Delta\overline{\ub}\ud M \\
   & = 0.
\end{align*}
It is easy to see that
\begin{multline*}
   \int_\Omega\left( (\ub\cdot\nabla)\ub + \left(\nabla\cdot\int_z^0\ub\ud\xi\right)
   \dfrac{\p\ub}{\p z}\right)\cdot(\Delta\ub)\ud\Omega \leq \\
   C\int_\Omega\left(|\ub||\nabla\ub| + \int_{-H}^0|\nabla\ub|\ud z|\p_z\ub|\right)
   |\Delta\ub|\ud\Omega.
\end{multline*}
Therefore we derive from \eqref{e79} that
\begin{multline}
   \dfrac{1}{2}\dfrac{\ud}{\ud t}|\nabla\ub|_{L^2(\Omega)}^2 +
   \dfrac{\mu}{2}|\Delta\ub|_{L^2(\Omega)}^2 +
   \nu|\nabla\p_z\ub|_{L^2(\Omega)}^2 \\
   + \mu_\delta|\Delta (I-P_{M,N})\ub|_{L^2(\Omega)}^2 +
   \nu_\delta|\nabla \p_z(I-P_{M,N})\ub|_{L^2(\Omega)}^2
   \leq \dfrac{1}{\mu}|\mathbf{F}|_{L^2}^2 + \dfrac{\mu}{4}|\Delta\ub|_{L^2}^2 \\
   + C\int_\Omega\left(|\ub||\nabla\ub| + \int_{-H}^0|\nabla\ub|\ud z|\p_z\ub|\right)
   |\Delta\ub|\ud\Omega.
   \label{e80}
\end{multline}
By the Young's inequality and the Sobolev interpolation inequality \eqref{a4} for $L^3$
functions, we have
\begin{align*}
   \int_\Omega |\ub||\nabla\ub||\Delta\ub|\ud\Omega
   &\leq C|\ub|_{L^6}|\nabla\ub|_{L^3}|\Delta\ub|_{L^2}\\
   &C\leq |\ub|_{L^6}|\nabla\ub|_{L^2}^{\frac{1}{2}}|\nabla\p_z\ub|_{L^2}^\frac{1}{2}
   |\Delta\ub|_{L^2} +
   C|\ub|_{L^6}|\nabla\ub|_{L^2}^{\frac{1}{2}}|\Delta\ub|_{L^2}^\frac{3}{2}.
\end{align*}
We appeal to the following inequality,
\begin{displaymath}
   \int_\Omega\left(\int_{-H}^0|\nabla\ub|\ud z\right)|f||g|\ud\Omega \leq
   C|f|_{L^2}|\ub|_{H^1}^{\frac{1}{2}}|\ub|_{H^2}^{\frac{1}{2}}
   |g|_{L^2}^{\frac{1}{2}}|g|_{H^1}^{\frac{1}{2}},
\end{displaymath}
to obtain (by setting $f=\Delta\ub$, $g=\p_z\ub$)
\begin{displaymath}
   \int_\Omega\left( \int_{-H}^0|\nabla\ub|\ud z\right)|\p_z\ub|
   |\Delta\ub|\ud\Omega \leq
   C|\nabla\ub|_{L^2}^{\frac{1}{2}}|\p_z\ub|_{L^2}^{\frac{1}{2}}
   |\nabla\p_z\ub|_{L^2}^{\frac{1}{2}}|\Delta\ub|_{L^2}^{\frac{3}{2}}.
\end{displaymath}
After these steps, we infer from \eqref{e80} that
\begin{multline}
   \dfrac{\ud}{\ud t}|\nabla\ub|_{L^2}^2 + {\mu}|\Delta\ub|_{L^2}^2
   + { }\\
   \nu|\nabla\p_z\ub|_{L^2}^2 + 2\mu_\delta|\Delta (I-P_{M,N})\ub|_{L^2}^2
   + 
   2\nu_\delta |\nabla \p_z(I-P_{M,N})\ub|_{L^2}^2 \leq \\
   \dfrac{4}{\mu}|\mathbf{F}|_{L^2}^2 + 
   C\left(|\ub|_{L^6}^4 + |\p_z\ub|_{L^2}^2|\nabla\p_z\ub|_{L^2}^2\right)
   |\nabla\ub|_{L^2}^2.
   \label{e81}
\end{multline}
We are ready to apply the Gronwall inequality, with use of \eqref{e77},
\eqref{e78}, to obtain
\begin{multline}
   |\nabla\ub(\cdot,t)|_{L^2}^2 + {\mu}\int_0^t|\Delta\ub|_{L^2}^2\ud s
   + \nu\int_0^t|\nabla\p_z\ub|_{L^2}^2\ud s + { }\\
   2\mu_\delta\int_0^t |\Delta (I-P_{M,N})\ub|_{L^2}^2\ud s 
   + 2\nu_\delta\int_0^t |\nabla \p_z(I-P_{M,N})\ub|_{L^2}^2\ud s \\
   \leq
   K_V(t),
   \label{e82}
\end{multline}
with
\begin{displaymath}
   K_V(t) =
   e^{C\left(K_6^{2} + K_2^2(t)\right)t + K_z^2(t)}\left(
   |\ub_0|_{H^1}^2 + \dfrac{4}{\mu}\int_0^t|\mathbf{F}(\cdot,s)|_{L^2}^2\ud s
   \right).
\end{displaymath}

Now that the key estimates are in place, with an argument
similar to that for Theorem \ref{t0} in Section \ref{s3},
we can show
\begin{theorem}\label{t1}
   For a given $T>0$, let
   $\mathbf{F}\in L^2(0,T;L^2(\Omega))$,
   $\ub_0\in V$. Then there exists a unique strong solution of the
   system \eqref{e1}--\eqref{e8} which depends continuously
   on the initial data.
\end{theorem}
The proof is omitted.

\subsection{Convergence of the solutions}\label{s2.3}
In this section we take $\mu_\delta \longrightarrow 0$, $\nu_\delta \longrightarrow 0$
in \eqref{e12a}, and study the convergence of the solution of the system.
We first rewrite \eqref{e12a} as follows:
\begin{multline}
   \dfrac{\p\ub^\delta}{\p t} + (\ub^\delta\cdot\nabla)\ub^\delta +
   (\nabla\cdot\int_z^0\ub^\delta\ud\xi)\dfrac{\p\ub^\delta}{\p z}
   +f\bk\times\ub^\delta + \dfrac{1}{\rho_0}\nabla p \\
   { } - \mu\Delta\ub^\delta
   - \nu\dfrac{\p^2\ub^\delta}{\p z^2} -
   \mu_\delta\Delta (I-P_{M,N})\ub^\delta -
   \nu_\delta\dfrac{\p^2}{\p z^2}(I-P_{M,N})\ub^\delta= \mathbf{F}.
   \label{e92}
\end{multline}
The superscript $\delta$ emphasizes the fact that the solution
$\ub^\delta$ depends on the spectral viscosity parameters $\mu_\delta$
and $\nu_\delta$, which themselves are determined by the grid resolution.
The proper relation between the parameters $\mu_\delta$ and $\nu_\delta$
and the grid resolution $\delta$ is the subject of a separate endeavor,
and will be presented elsewhere. In this work we focus on the behavior
of the solution $\ub_\delta$ of system \eqref{e92} as $\mu_\delta$
and $\nu_\delta$ tend to 0 (corresponding to the scenario when
the grid becomes finer and finer).

Let $\ub$ be the solution of the primitive equations without artificial
viscosities, that is, $\ub$ is the solution of
\begin{multline}
   \dfrac{\p\ub}{\p t} + (\ub\cdot\nabla)\ub +
   (\nabla\cdot\int_z^0\ub\ud\xi)\dfrac{\p\ub}{\p z}
   +f\bk\times\ub + { }\\
   \dfrac{1}{\rho_0}\nabla p-
   \mu\Delta\ub
   - \nu\dfrac{\p^2\ub}{\p z^2} = \mathbf{F}.
   \label{e93}
\end{multline}
It can be shown  that the system \eqref{e93} plus \eqref{e4}--\eqref{e8}
has unique global strong solutions, provided that the initial
data and the forcing are sufficiently smooth. See \cite{CT07, Ko07, KuZi07}
(These works use boundary conditions different from ours, but the case with
periodic boundary
conditions can be handled as well.)

To study the convergence of the solution of system \eqref{e92},
we subtract \eqref{e93} from \eqref{e92}, let $\vb^\delta =
\ub^\delta - \ub$, and we have
\begin{multline}
   \dfrac{\p\vb^\delta}{\p t} + (\ub^\delta\cdot\nabla)\vb^\delta +
   (\nabla\cdot\int_z^0\ub^\delta\ud\xi)\dfrac{\p\vb^\delta}{\p z}
   + (\vb^\delta\cdot\nabla)\ub +
   (\nabla\cdot\int_z^0\vb^\delta\ud\xi)\dfrac{\p\ub}{\p z}
   + { }\\
   f\bk\times\vb^\delta
    - \mu\Delta\vb^\delta
   - \nu\dfrac{\p^2\vb^\delta}{\p z^2} - \mu_\delta\Delta (I-P_{M,N})\ub^\delta -
   \nu_\delta\dfrac{\p^2}{\p z^2}(I-P_{M,N})\ub^\delta= 0.
   \label{e94}
\end{multline}
We take the inner product of \eqref{e94} with $\vb^\delta$, and
integrate by parts over $\Omega$ to obtain
\begin{multline}
   \dfrac{1}{2}\dfrac{\ud}{\ud t}|\vb^\delta|_{L^2}^2 +
   \int_{\Omega}\left( (\ub^\delta\cdot\nabla)\vb^\delta +
    (\nabla\cdot\int_z^0\ub^\delta\ud\xi)\dfrac{\p\vb^\delta}{\p z}\right)\cdot
    \vb^\delta\ud\Omega + { } \\
   \int_{\Omega}\left( (\vb^\delta\cdot\nabla)\ub+
    (\nabla\cdot\int_z^0\vb^\delta\ud\xi)\dfrac{\p\ub}{\p z}\right)\cdot
    \vb^\delta\ud\Omega + { }\\
    \int_\Omega f\bk\times\vb^\delta\cdot\vb^\delta\ud\Omega +
   \mu|\nabla\vb^\delta|^2 + \nu\left|\dfrac{\p \vb^\delta}{\p z}\right|^2 \\
   =
   -\mu_\delta(\nabla (I-P_{M,N})\ub^\delta,\,\nabla\vb^\delta)
   -\nu_\delta\left(\dfrac{\p}{\p z}(I-P_{M,N})\ub,\,\dfrac{\p}{\p z}\vb^\delta\right).
   \label{e95}
\end{multline}
We note that
\begin{equation}
   \int_\Omega f\bk\times\vb^\delta\cdot\vb^\delta\ud\Omega = 0.
   \label{e96}
\end{equation}
We can also verify by integration by parts, and using the boundary conditions
\eqref{e4}--\eqref{e6} for $\ub$, $\ub^\delta$ (and therefore for $\vb^\delta$),
that
\begin{equation}
   \int_{\Omega}\left( (\ub^\delta\cdot\nabla)\vb^\delta +
    (\nabla\cdot\int_z^0\ub^\delta\ud\xi)\dfrac{\p\vb^\delta}{\p z}\right)\cdot
    \vb^\delta\ud\Omega = 0.
    \label{e97}
\end{equation}
By Holder's inequality and the interpolation inequalities \eqref{a4} and \eqref{a5} in $R^3$, we find that
\begin{align}
   &\left|\int_{\Omega}\left( (\vb^\delta\cdot\nabla)\ub\right)\cdot\vb^\delta\ud\Omega \right|
   \leq C|\nabla\ub|_{L^2}|\vb^\delta|_{L^3}|\vb^\delta|_{L^6},\nonumber\\
   &\left|\int_{\Omega}\left( (\vb^\delta\cdot\nabla)\ub\right)\cdot\vb^\delta\ud\Omega \right|
   \leq C|\nabla\ub|_{L^2}|\vb^\delta|_{L^2}^{\frac{1}{2}}
   |\nabla\vb^\delta|_{L^2}^{\frac{3}{2}}.\label{e98}
\end{align}
Using Holder's inequality again we find that
\begin{align*}
   &\phantom{\leq}\left|\int_\Omega \left(\nabla\cdot\int_z^0\vb^\delta\ud\xi\right)
   \dfrac{\p\ub}{\p z}\cdot\vb^\delta\ud\Omega\right|\\
   &\leq
   C\int_M\left(\int_{-H}^0|\nabla\vb^\delta|\ud z\right)
   \left(\int_{-H}^0\left|\dfrac{\p\ub}{\p z}\right|
   |\vb^\delta|\ud z\right)\ud x\ud y  \\
   &\leq
   C\int_M\left(\int_{-H}^0|\nabla\vb^\delta|\ud z\right)
   \left(\int_{-H}^0\left|\dfrac{\p\ub}{\p z}\right|^2\ud z\right)^
   {\frac{1}{2}}\left(\int_{-H}^0
   |\vb^\delta|^2\ud z\right)^{\frac{1}{2}}\ud x\ud y  \\
   &\leq
   C\left(\int_M\left(\int_{-H}^0|\nabla\vb^\delta|\ud z\right)^2
   \ud x\ud y\right)^{\frac{1}{2}}
   \left(\int_M
   \left(\int_{-H}^0\left|\dfrac{\p\ub}{\p z}\right|^2\ud z\right)^2
   \ud x\ud y\right)^{\frac{1}{4}}\times { }\\
   &\phantom{\leq} \left(\int_M\left(\int_{-H}^0
   |\vb^\delta|^2\ud z\right)^2\ud x\ud y\right)^{\frac{1}{4}}.
\end{align*}
We first note that
\begin{align*}
   \left(\int_M\left(\int_{-H}^0|\nabla\vb^\delta|\ud z\right)^2
   \ud x\ud y\right)^{\frac{1}{2}} &\leq
   C\left(\int_M\int_{-H}^0|\nabla\vb^\delta|^2\ud z
   \ud x\ud y\right)^{\frac{1}{2}} \\
   &= C|\nabla\vb^\delta|_{L^2(\Omega)}
\end{align*}
By the integral Minkowski inequality \eqref{a6} and the interpolation inequality \eqref{a2}
in $R^2$, we find that
\begin{align*}
   \left(\int_M
   \left(\int_{-H}^0\left|\dfrac{\p\ub}{\p z}\right|^2\ud z\right)^2
   \ud x\ud y\right)^{\frac{1}{2}} &\leq
   C\int_{-H}^0
   \left(\int_M\left|\dfrac{\p\ub}{\p z}\right|^4\ud x\ud y\right)^
   {\frac{1}{2}}\ud z\\
   &\leq C\int_{-H}^0\left|\dfrac{\p\ub(\cdot,z)}{\p z}\right|_{L^4(M)}^2\ud z\\
   &\leq C\int_{-H}^0 \left|\dfrac{\p\ub(\cdot,z)}{\p z}\right|_{L^2(M)}
   \left|\nabla\dfrac{\p\ub(\cdot,z)}{\p z}\right|_{L^2(M)}\ud z\\
   &\leq C\left|\dfrac{\p\ub}{\p z}\right|_{L^2(\Omega)}
   \left|\nabla\dfrac{\p\ub}{\p z}\right|_{L^2(\Omega)}.
\end{align*}
Similarly, we have
\begin{equation*}
\left(\int_M\left(\int_{-H}^0
   |\vb^\delta|^2\ud z\right)^2\ud x\ud y\right)^{\frac{1}{2}}
   \leq C|\vb^\delta|_{L^2}|\nabla\vb^\delta|_{L^2}.
\end{equation*}
Therefore,
\begin{equation}
   \left|\int_\Omega \left(\nabla\cdot\int_z^0\vb^\delta\ud\xi\right)
   \dfrac{\p\ub}{\p z}\cdot\vb^\delta\ud\Omega\right|
   \leq C\left|\dfrac{\p\ub}{\p z}\right|_{L^2(\Omega)}^{\frac{1}{2}}
   \left|\nabla\dfrac{\p\ub}{\p z}\right|_{L^2(\Omega)}^{\frac{1}{2}}
   |\vb^\delta|_{L^2}^{\frac{1}{2}}|\nabla\vb^\delta|_{L^2}^{\frac{3}{2}}.
   \label{e99}
\end{equation}
With \eqref{e96}--\eqref{e99}, we derive from \eqref{e95} that
\begin{multline}
   \dfrac{1}{2}\dfrac{\ud}{\ud t}|\vb^\delta|_{L^2}^2 +
    \mu|\nabla\vb^\delta|^2 + \nu\left|\dfrac{\p \vb^\delta}{\p z}\right|^2
    \leq
   \mu_\delta\left|\nabla (I-P_{M,N})\ub^\delta\right|\left|\nabla\vb^\delta\right| +
   \\
   \nu_\delta\left|\dfrac{\p}{\p z}(I-P_{M,N})\ub\right|\left|
   \dfrac{\p}{\p z}\vb^\delta\right| +
   C\left(\left|\dfrac{\p\ub}{\p z}\right|_{L^2(\Omega)}^{\frac{1}{2}}
   \left|\nabla\dfrac{\p\ub}{\p z}\right|_{L^2(\Omega)}^{\frac{1}{2}}
   \right)
   |\vb^\delta|_{L^2}^{\frac{1}{2}}|\nabla\vb^\delta|_{L^2}^{\frac{3}{2}} +{ }\\
   C|\nabla\ub|_{L^2}
   |\vb^\delta|_{L^2}^{\frac{1}{2}}|\nabla\vb^\delta|_{L^2}^{\frac{3}{2}}.
   \label{e100}
\end{multline}
By Young's inequality, we have
\begin{multline}
   \dfrac{\ud}{\ud t}|\vb^\delta|_{L^2}^2 +
    \mu|\nabla\vb^\delta|^2 + \nu\left|\dfrac{\p \vb^\delta}{\p z}\right|^2
    \leq
    2\dfrac{\mu_\delta^2}{\nu}\left|\nabla (I-P_{M,N})\ub^\delta\right|^2 +
   \\
   \dfrac{\nu_\delta^2}{\nu}\left|\dfrac{\p}{\p z}(I-P_{M,N})\ub\right|^2 +
   C\left(\left|\dfrac{\p\ub}{\p z}\right|_{L^2(\Omega)}^{2}
   \left|\nabla\dfrac{\p\ub}{\p z}\right|_{L^2(\Omega)}^{2} +
   |\nabla\ub|_{L^2}^4
   \right)
   |\vb^\delta|_{L^2}^2.
   \label{e101}
\end{multline}
We notice that the a priori estimates obtained in the previous section
are independent of $\mu_\delta$ and $\nu_\delta$. We apply the Gronwall
inequality to \eqref{e101}, and, utilizing the a priori estimates
\eqref{e37}, \eqref{e78} and \eqref{e82}, we obtain
\begin{multline}
   |\vb^\delta(\cdot,t)|_{L^2}^2 +
    \int_0^t\mu|\nabla\vb^\delta(\cdot,s)|^2\ud s +
    \nu\int_0^t\left|\dfrac{\p \vb^\delta}{\p z}(\cdot,s)\right|^2\ud s \\
    \leq e^{K_z^2(t)+K_V^2(t)}\left(\dfrac{2\mu_\delta}{\mu} +
    \dfrac{\nu_\delta}{\nu}\right)K_1(t).
    \label{e102}
 \end{multline}

 We have just proved
 \begin{theorem}\label{t2}
    Let $T>0$ be given, and the other assumptions be the same as in
    Theorem \ref{t1}. Then, as $\mu_\delta\,\&\,\nu_\delta \longrightarrow 0$,
    \begin{align}
       &\ub^\delta - \ub \sim \sqrt{\mu_\delta + \nu_\delta} \longrightarrow 0
       \quad\textrm{ in }L^\infty(0,T;L^2(\Omega)),\\
       &\ub^\delta - \ub \sim \sqrt{\mu_\delta + \nu_\delta} \longrightarrow 0
       \quad\textrm{ in }L^2(0,T;H^1(\Omega)).
    \end{align}
 \end{theorem}

\appendix
\section{Some inequalities}
We list here some functional inequalities that are frequently
used in this paper. \\
\noindent{\it $L^p$ interpolation inequality}\\
Let $\Omega \subset \mathbb{R}^3$, and  $1\leq p_1\leq p\leq p_2$, $p_1 \neq p_2$.
Let $u\in L^{p_1}(\Omega)\cap L^{p_2}(\Omega)$. Then $u\in L^p(\Omega)$, and
\begin{equation}
    ||u||_{L^p} \leq ||u||_{L^{p_1}}^{s_1} ||u||_{L^{p_2}}^{s_2},
    \label{a1}
\end{equation}
with
\begin{displaymath}
    s_1 = \dfrac{p_1}{p}\dfrac{p_2-p}{p_2-p_1},\qquad
    s_2 = \dfrac{p_2}{p}\dfrac{p-p_1}{p_2-p_1}.
\end{displaymath}

\vspace{5mm}
\noindent{\it Ladyzhenskaya/Sobolev inequalities in $\mathbb{R}^2$}\\
Let $M\subset\mathbb{R}^2$ be a bounded domain with piecewise smooth
boundaries. For each $\phi\in H^1(M)$, the following inequalities hold:
\begin{align}
    &||\phi||_{L^4(\Omega)} \leq C||\phi||_{L^2(\Omega)}^\frac{1}{2}
    ||\phi||_{H^1(\Omega)}^\frac{1}{2},\label{a2}\\
    &||\phi||_{L^8(\Omega)} \leq C||\phi||_{L^6(\Omega)}^\frac{3}{4}
    ||\phi||_{H^1(\Omega)}^\frac{1}{4}.\label{a3}
\end{align}

\vspace{5mm}
\noindent{\it Ladyzhenskaya/Sobolev inequalities in $\mathbb{R}^3$}\\
Let $\Omega\subset\mathbb{R}^3$ be a bounded domain with piecewise smooth
boundaries. For each $\phi\in H^1(\Omega)$, the following inequalities hold:
\begin{align}
    &||\phi||_{L^3(\Omega)} \leq C||\phi||_{L^2(\Omega)}^\frac{1}{2}
    ||\phi||_{H^1(\Omega)}^\frac{1}{2},\label{a4}\\
    &||\phi||_{L^6(\Omega)} \leq C||\phi||_{H^1(\Omega)}.\label{a5}
\end{align}

\vspace{5mm}
\noindent{\it Minkowski integral inequality (for $p \geq 1$)}\\
Let $\Omega_1\subset \mathbb{R}^{m_1}$ and $\Omega_2\subset \mathbb{R}^{m_2}$
be two measurable sets, with $m_1$ and $m_2$ being positive integers.
Let $f(\xi,\eta)$ be a measurable function over $\Omega_1\times\Omega_2$.
Then
\begin{equation}
    \left(\int_{\Omega_1}\left(\int_{\Omega_2}
    |f(\xi,\eta)|\ud\eta\right)^p\ud\xi\right)^{\frac{1}{p}} \leq
    \int_{\Omega_2}\left(\int_{\Omega_1}|f(\xi,\eta)|^p
    \ud\xi\right)^{\frac{1}{p}}\ud\eta.
    \label{a6}
\end{equation}
The inequality \eqref{a1} can be verified by the Holder's inequality. For
\eqref{a2}--\eqref{a5} we
refer to such classical texts as \cite{Ad75, Te01}. A proof of \eqref{a6}
can be found in \cite{HLP88}.

\section*{Acknowledgment}
Q. Chen and M. Gunzburger are supported by the US Department of 
Energy grant number DE-SC0002624 as part of the 
{\em Climate Modeling: Simulating Climate at Regional Scale} program.
This work is supported in part by grants from the National Science
Foundation. Wang acknowledges helpful conversation with Ning Ju.

\bibliographystyle{amsplain}
\bibliography{references}

\end{document}